\documentclass[arxiv,reqno,bibliography=totoc,twoside,a4paper,12pt]{amsart}

\usepackage[utf8]{inputenc}
\usepackage{amsmath,bbm}
\usepackage{amsfonts}
\usepackage{amssymb}
\usepackage{amsthm}
\usepackage{float}
\usepackage{mathtools}
\usepackage{csquotes}
\usepackage{upgreek,esint,enumerate}
\usepackage[numbers]{natbib}
\usepackage{dsfont,ednotes}
\usepackage{enumitem}
\usepackage{graphicx}

\usepackage[colorlinks=true,linkcolor=blue,citecolor=blue, urlcolor=blue]{hyperref}
\usepackage{cleveref}

\usepackage[scaled]{helvet} 
\usepackage{courier} 
\usepackage[mathbf]{euler}
\usepackage[dvipsnames]{xcolor}

\usepackage[driver=pdftex,margin=3cm,heightrounded=true,centering]{geometry}

\usepackage{tikz-cd}
\usetikzlibrary{arrows,matrix,shapes,decorations.text, mindmap,shapes.misc,decorations.markings}

\usepackage{metalogo}
\usepackage{commands} 
\usepackage[normalem]{ulem}

\theoremstyle{plain}
\newtheorem{theo}{Theorem}[section]

\newtheorem{lem}[theo]{Lemma}
\newtheorem{cor}[theo]{Corollary}
\newtheorem{prop}[theo]{Proposition}
\newtheorem{defi}[theo]{Definition}
\newtheorem{rem}[theo]{Remark}

\usepackage[toc,page]{appendix}

\numberwithin{equation}{section}

\begin{document}

\title[$L^2$-Gamma index theorem for spacetimes]
{$L^2$-Gamma index theorem for spacetimes}

\author{Orville Damaschke}
\address{Universit\"at Oldenburg,
26129 Oldenburg,
Germany}

\author{Boris Vertman}
\address{Universit\"at Oldenburg,
26129 Oldenburg,
Germany}
\email{boris.vertman@uni-oldenburg.de}

\subjclass[2000]{58J20; 58J30, 58J28}
\date{\today}

\begin{abstract}
{We establish an $L^2$-Gamma index theorem for the Dirac operator on a globally hyperbolic manifold $M$ 
with Cauchy hypersurface $\Sigma$ being a Galois covering of a compact smooth manifold
with Galois group $\Gamma$. Our argument rewrites the $L^2$-Gamma index
in terms of the spectral flow, which is then connected to the usual geometric expressions. This extends the 
work of B\"ar and Strohmaier to some non-compact Cauchy hypersurfaces. The analysis here is 
based on intermediate results by the first author on $L^2$-Gamma Fredholm properties of 
the Dirac operator in \cite{OD1} and \cite{OD2}.}
\end{abstract}

\maketitle
\tableofcontents

\section{Introduction and statement of the main result}\label{chap:Intro}

\subsection{Geometric setting}\label{setting-subsection} 

We recall briefly the notations and the setting from \cite{OD1} and \cite{OD2}, especially \cite[Sec.5]{OD2}.
We consider a globally hyperbolic temporally compact manifold 
$$
M = [t_1,t_2] \times \Sigma,
$$
of even dimension $(n + 1)$, i.e. $n = \dim \Sigma$ is odd. The spacelike Cauchy hypersurface 
$\Sigma$ is assumed to be a Galois covering with respect to a discrete Galois group $\Gamma$ 
acting by deck transformations with compact smooth quotient $\Sigma_\Gamma = \Sigma / \Gamma$. 
Thus $M$ becomes a Galois covering as well, with Galois group $\Gamma$ and globally 
hyperbolic base $M_\Gamma = [t_1,t_2] \times \Sigma_\Gamma$.
The spacetime $M$ is equipped with a Lorentzian metric 
$$
\met = -N^2 dt^2 + \met_t,
$$
where $N \in C^\infty(M,\R)$ is nowhere vanishing, and $\met_t, t \in [t_1,t_2]$, is a smooth family
of Riemannian $\Gamma$-invariant metrics on the Galois covering $\Sigma$.
We will denote the \emph{flipped} Riemannian metric on $M$ by $\check{\met} = N^2 dt^2 + \met_t$.

\subsection{The spin and spin$^c$ Dirac operator} 

Assume that $M$ is spin and consider the spinor bundle $\mathscr{S}(M)$.
We may also assume more generally that $M$ is only spin$^c$. In that case 
$\mathscr{S}(M)$ is only locally defined. The spin$^c$ structure defines a Hermitian 
line bundle $L$ with a locally defined square root bundle $L^{\frac{1}{2}}$, such that
the twisted spin$^c$ spinor bundle $\mathscr{S}_L(M) = \mathscr{S}(M) \otimes L^{\frac{1}{2}}$
is well-defined globally.\medskip

Since $\dim M$ is even, $\mathscr{S}_L(M)$ decomposes into a direct sum of 
half spinor bundles $\mathscr{S}_L(M) = \mathscr{S}_L^+(M) \oplus \mathscr{S}_L^-(M)$. The spinor bundles
are $\Gamma$-equivariant and the spin$^c$ Dirac operators acting on sections of spinor bundles
$$
D_\pm: C^\infty(M, \mathscr{S}_L^\pm(M)) \to C^\infty(M, \mathscr{S}_L^\mp(M)),
$$ 
are $\Gamma$-equivariant lifts of the corresponding spin$^c$ Dirac operators on the 
base $M_\Gamma$, which is also spin$^c$. We also consider a Hermitian $\Gamma$-equivariant
vector bundle $E$ over $M$. We set $\mathscr{S}_{L,E}^\pm(M) := \mathscr{S}_{L}^\pm(M) \otimes E$ and 
consider the twisted spin$^c$ Dirac operators acting on sections of the half spinor bundles $\mathscr{S}_{L,E}^\pm(M)$
$$
D^{E_L}_\pm: C^\infty(M, \mathscr{S}_{L,E}^\pm(M)) \to C^\infty(M, \mathscr{S}_{L,E}^\mp(M)).
$$ 
The restrictions of $\mathscr{S}_{L,E}^\pm(M)$ to $\Sigma_t := \{t\} \times \Sigma$
can be naturally identified with the twisted spin$^c$ spinor bundle $\mathscr{S}_{L,E}(\Sigma_t)$ 
of $\Sigma_t$. Let $\nu$ be the past-directed timelike vector field on $(M,g)$ with $g(\nu, \nu) = -1$,  	
perpendicular to all $\Sigma_t$. We write $\beta$ for the Clifford multiplication on $\mathscr{S}_L(M)$ by $\nu$.
Then the twisted spin$^c$ Dirac operators along each $\Sigma_t$ take the following form			
\begin{align}\label{D-product-boundary}
D^{E_L}_\pm = - (\beta \otimes \mathds{1}_{E_L}) \left(\nabla^{\mathscr{S}_{L,E}(\Sigma_t)}_\nu 
\pm i A^{E_L}_t - \frac{n  H_t }{2}\right),
\end{align}
where $H_t$ is the mean curvature of $\Sigma_t$ with respect to $\nu$,
and $A^{E_L}_t$ is the twisted Dirac operator on $\Sigma_t$. See 
\cite[\S 3.2]{OD1} for further details.

\subsection{Generalized APS and aAPS boundary conditions} 

We shall impose boundary conditions on the two boundary components of $M$,
$\Sigma_1 \equiv \Sigma_{t_1}$ and $\Sigma_2 \equiv \Sigma_{t_2}$.
We shall also abbreviate $A^{E_L}_j \equiv A^{E_L}_{t_j}, j=1,2$.
The tangential operators $A^{E_L}_1, A^{E_L}_2$ act on sections of 
$\mathscr{S}_{L,E}(\Sigma)$. They are essentially self-adjoint with real spectrum 
and associated spectral projections $P_I(t_j):= \Chi_I(A^{E_L}_j)$ for any interval $I \subset \R$ and $j=1,2$,  see \cite[Sec. 7.1]{OD2} for details. 	\medskip
											
We denote the space of square-integrable sections of $\mathscr{S}_{L,E}(\Sigma_t)$ with respect to a $\Gamma$-invariant measure 
on the Galois covering $\Sigma_t$ (e.g. induced by the $\Gamma$-invariant metric $\met_t$) by $L^2(\Sigma_t,\mathscr{S}_{L,E}(\Sigma_t))$ and define 
$L^2_I(\Sigma_t,\mathscr{S}_{L,E}(\Sigma_t))$ to be the range of $L^2(\Sigma_t,\mathscr{S}_{L,E}(\Sigma_t))$ under $P_I(t)$,  where we also use the notation
$\range{P_I(t_j)}$.  We also write $L^2(M,\mathscr{S}^{\pm}_{L,E}(M))$ for square-integrable $\Gamma$-sections over the Galois 
covering $M$. These spaces are examples of free $\Gamma$-Hilbert modules.

\begin{defi}
The generalized APS and aAPS boundary conditions
 are defined for any given $a_1, a_2 \in \R$ by the following maps 								
\begin{equation*}
\begin{split}
\mathbb{P}_{+} := \, &(P_{\intervallro{a_1}{\infty}{}}(t_1)\circ\rest{\Sigma_1})\oplus (P_{\intervalllo{-\infty}{a_2}{}}(t_2)\circ\rest{\Sigma_{2}}): \\ 
&C^0([t_1,t_2],  L^2(\Sigma_\bullet, \mathscr{S}_{L,E}(\Sigma_\bullet))) \to L^2_{[a_1,\infty)}(\Sigma_1, \mathscr{S}_{L,E}(\Sigma_1)) 
\oplus L^2_{(-\infty,a_2]}(\Sigma_2, \mathscr{S}_{L,E}(\Sigma_2)), \\
\mathbb{P}_{-} := \, &(P_{\intervallo{-\infty}{a_1}{}}(t_1)\circ\rest{\Sigma_1})\oplus (P_{\intervallo{a_2}{\infty}{}}(t_2)\circ\rest{\Sigma_{2}}): \\ 
&C^0([t_1,t_2],  L^2(\Sigma_\bullet, \mathscr{S}_{L,E}(\Sigma_\bullet))) \to L^2_{(-\infty,a_1)}(\Sigma_1, \mathscr{S}_{L,E}(\Sigma_1)) 
\oplus L^2_{(a_2,\infty)}(\Sigma_2, \mathscr{S}_{L,E}(\Sigma_2)).
\end{split}
\end{equation*} 
Here, $C^0([t_1,t_2],  L^2(\Sigma, \mathscr{S}_{L,E}(\Sigma)))$ is the space of continuous functions
in $[t_1,t_2]$ with values in the space of square-integrable sections of the spinor bundle $\mathscr{S}_{L,E}(\Sigma)$
over $\Sigma$. For any $j=1,2$, the notation $\rest{\Sigma_j}$ indicates the restriction of such functions to 
$t_j$. 
\end{defi}

We shall use spaces of finite energy $\Gamma$-sections of order zero, 
introduced in \cite{OD2}:
\begin{equation*}
\begin{split}
&FE^0_{\Gamma}(M,D_{\pm}^{E_L}) := 
\SET{u \in C^0([t_1,t_2],  L^2(\Sigma_\bullet, \mathscr{S}_{L,E}(\Sigma_\bullet))): D_{\pm}^{E_L}u \in L^2(M, \mathscr{S}_{L,E}^\mp(M))}, \\ 
&FE^0_{\Gamma,\mathrm{APS}(a_1,a_2)}(M,D_{\pm}^{E_L}) :=
\SET{u \in FE^0_{\Gamma}(M,D_{\pm}^{E_L}): \mathbb{P}_{\pm} u = 0},\\ 
&FE^0_{\Gamma,\mathrm{aAPS}(a_1,a_2)}(M,D_{\pm}^{E_L}) :=
\SET{u \in FE^0_{\Gamma}(M,D_{\pm}^{E_L}): \mathbb{P}_{\mp} u = 0}.
\end{split}
\end{equation*}

The Dirac operators with generalized APS and aAPS boundary conditions are
\begin{equation}\label{gammafredoperres}
\begin{split}
D^{E_L}_{\pm,\mathrm{APS}(a_1,a_2)}\,&:\,FE^0_{\Gamma,\mathrm{APS}(a_1,a_2)}(M,D_{\pm}^{E_L})
\,\rightarrow L^2(M, \mathscr{S}^{\mp}_{L,E}(M)),\\
D^{E_L}_{\pm,\mathrm{aAPS}(a_1,a_2)}\,&:\,FE^0_{\Gamma,\mathrm{aAPS}(a_1,a_2)}(M,D_{\pm}^{E_L})
\,\rightarrow L^2(M, \mathscr{S}^{\mp}_{L,E}(M)).
\end{split}
\end{equation}

The operators in \eqref{gammafredoperres} are $\Gamma$-Fredholm 
by \cite[Theorem 9.1]{OD2}. Their $\Gamma$-indices are given by 
$\Gamma$-indices of
some wave evolution operators, reviewed below in \S \ref{chap:recap} and
are related among each other as follows.

\begin{prop}
The operators \eqref{gammafredoperres} are 
$\Gamma$-Fredholm with $\Gamma$-indices
\begin{align*}
\Index_{\Gamma}(D^{E_L}_{+,\mathrm{APS}(a_1,a_2)})&= 
\Index_{\Gamma}(D^{E_L}_{-,\mathrm{aAPS}(a_1,a_2)})= -\Index_{\Gamma}(D^{E_L}_{+,\mathrm{aAPS}(a_1,a_2)})\\
&=-\Index_{\Gamma}(D^{E_L}_{-,\mathrm{APS}(a_1,a_2)}).
\end{align*}
Moreover,  denoting with $\dim_\Gamma$ the $\Gamma$-dimension of a $\Gamma$-Hilbert (sub-)module,  the $\Gamma$-indices for the generalized $\mathrm{APS}(a_1,a_2)$ 
and the (classical) APS boundary conditions (i.e. $\mathrm{APS}(0,0)$) are related by
\begin{align*}
\Index_{\Gamma}(D^{E_L}_{+,\mathrm{APS}(a_1,a_2)}) &= 
\Index_{\Gamma}(D^{E_L}_{+,\mathrm{APS}(0,0)}) \\ &+ \Chi_{\SET{a_2 < 0}}
\dim_\Gamma\left(\range{P_{(a_2,0]}(t_2)}\right)-\Chi_{\SET{a_2 > 0}}
\dim_\Gamma\left(\range{P_{(0,a_2]}(t_2)}\right)\\
&+ \Chi_{\SET{a_1 > 0}}\dim_\Gamma\left(\range{P_{[0,a_1)}(t_1)}\right)-\Chi_{\SET{a_1 < 0}}
\dim_\Gamma\left(\range{P_{[a_1,0)}(t_1)}\right),
\end{align*}
where we abbreviated $\Chi_{\SET{a>0}} := \Chi_{(0,\infty)}(a)$ and $\Chi_{\SET{a<0}}:= \Chi_{(-\infty,0)}(a)$.
\end{prop}

The proof of these statements follows from Theorem \ref{qgengammafred}
(Theorem 9.2 in \cite{OD2}) and Theorem \ref{indexDgaAPStwist} (Theorem 9.3 in \cite{OD2}).

\subsection{Statement of the main result} 

The $\Gamma$-equivariant Hermitian vector bundles 
$E$ and $L$ over the Galois $\Gamma$-covering $M$ 
are lifts of Hermitian vector bundles $E_\Gamma$ and $L_\Gamma$ over the 
base $M_\Gamma$. Similarly $\met$ is the lift of the Lorentzian metric 
$\met_\Gamma$ on $M_\Gamma$. The operators $A^{E_L}_j$ are lifts
of the corresponding operators $\underline{A}^{E_L}_j$ on the compact quotients 
$\Sigma_j / \Gamma$, $j=1,2$. Our main result here is an explicit index formula, see
\S \ref{chap:geoind} below for details.

\begin{theo}\label{mainres2}
The index of $D^{E_L}_{+,\mathrm{APS}} \equiv D^{E_L}_{+,\mathrm{APS}(0,0)}$
is given explicitly as follows
\begin{equation*}\begin{split}
\Index_{\Gamma}(D^{E_L}_{+,\mathrm{APS}}) 
&= \int_{M_\Gamma}\hat{\mathbbm{A}}\left(M_\Gamma\right)\wedge\expe{\mathrm{c}_1\left(L_\Gamma\right)/2}\wedge \mathrm{ch}\left(E_\Gamma\right)
\\ &+\int_{\bound M_\Gamma} T\widehat{\mathbbm{A}}(g)\wedge\expe{\mathrm{c}_1\left(L_\Gamma\vert_{\bound M}\right)/2}\wedge \mathrm{ch}\left(E_\Gamma\vert_{\bound M}\right))\\
&-\xi_\Gamma(A_{1}^{E_L})+\xi_\Gamma(A_{2}^{E_L}) - \dim_\Gamma \ker A^{E_L}_{t_2}.
\end{split}\end{equation*}
Here, we have used the following notation: \medskip

\begin{centering}
\begin{enumerate} 
\item \quad $\bound M_\Gamma := \Sigma_1/\Gamma\sqcup \Sigma_2/\Gamma$, \smallskip

\item \quad $\hat{\mathbbm{A}}\left(M_\Gamma\right)$ is the $\check{A}$-genus of $(M_\Gamma,\met_\Gamma)$, \smallskip

\item \quad $\mathrm{ch}\left(E_\Gamma\right)$ is the Chern character of $E_\Gamma$, \smallskip

\item \quad $\mathrm{c}_1\left(L_\Gamma\right)$ is the first Chern class of the line bundle $L_\Gamma$, \smallskip

\item \quad $T\widehat{\mathbbm{A}}(g)$ is pullback of the transgression form to $\partial M$, \smallskip

\item \quad $\xi_\Gamma$ is the reduced $\Gamma$-eta invariant, \smallskip

\item \quad $\dim_\Gamma\ker A^{E_L}_{t_2}$ is the $\Gamma$-dimension of the kernel of $A^{E_L}_{t_2}$ as $\Gamma$-Hilbert submodule.
\end{enumerate}
\end{centering}

\end{theo}

We prove Theorem \ref{mainres2} in the same way as it has been done in 
\cite{BaerStroh} for compact Cauchy hypersurfaces. The $\Gamma$-indices 
turn out to be the spectral flow of a family of Dirac operators on the hypersurfaces $\Sigma_t$. This provides algebraic $\Gamma$-index expressions for each Dirac operator \clef{gammafredoperres}. Based on known results from the Riemannian setting, we first express the geometric $\Gamma$-index of a suitable chosen Riemannian Dirac operator in terms of the $\Gamma$-index for the same family of hypersurface Dirac operators. This allows to connect the algebraic Lorentzian $\Gamma$-index with geometry. It is then left to show that the geometric expression, rewritten with the Lorentzian connection, does not depend on the chosen auxiliary Riemannian setting such that the resulting geometric $\Gamma$-index formula is purely given by those data, only required from the Lorentzian setting.

It has been brought to our attention by Alexander Strohmaier,  that the recent local index theorem in \cite{BaerStroh2} provides an alternative ansatz to our result.

\section{Review of previous results on Gamma-Fredholmness}\label{chap:recap}

Parallel to the definition of the finite energy sections $FE^0_{\Gamma}(M,D_{\pm}^{E_L})$
above, we define the space of finite energy kernel solutions
\begin{equation*}
FE^0_{\Gamma}(M,\ker D_{\pm}^{E_L}) := 
\SET{u \in C^0([t_1,t_2],  L^2(\Sigma, \mathscr{S}_{L,E}(\Sigma))): D_{\pm}^{E_L}u =0}
\end{equation*}
Then, by a crucial result by the first author,  \cite[Corollary 6.2]{OD2}, 
the restriction of functions in $C^0([t_1,t_2],  L^2(\Sigma_\bullet, \mathscr{S}_{L,E}(\Sigma_\bullet)))$
to $\Sigma_t = \{t\}\times \Sigma$ defines a $\Gamma$-isomorphism of $\Gamma$-Hilbert modules
\begin{equation*}
\rest{t}^\pm: FE^0_{\Gamma}(M,\ker D_{\pm}^{E_L}) \to L^2(\Sigma_t, \mathscr{S}_{L,E}(\Sigma_t)).
\end{equation*}
The notion of $\Gamma$-isomorphisms and $\Gamma$-Hilbert modules
are reviewed in e.g.  \cite[Def. 3.2 and Prop. 3.3]{OD2}.  This allows us to 
define two central $\Gamma$-isomorphisms of $\Gamma$-Hilbert modules,
cf.  \cite[Def.6.3]{OD2}
\begin{equation}\label{Q-def} 
\begin{split}
&Q(t_2,t_1):= \rest{t_2}^+ \circ (\rest{t_1}^+)^{-1}: L^2(\Sigma_1, \mathscr{S}_{L,E}(\Sigma_1))
\to L^2(\Sigma_2, \mathscr{S}_{L,E}(\Sigma_2)), \\
&\tilde{Q}(t_2,t_1):= \rest{t_2}^- \circ (\rest{t_1}^-)^{-1}: L^2(\Sigma_1, \mathscr{S}_{L,E}(\Sigma_1))
\to L^2(\Sigma_2, \mathscr{S}_{L,E}(\Sigma_2)).
\end{split}
\end{equation}
Using spectral projections $P_I(t_j):= \Chi_I(A^{E_L}_j)$ we decompose,  as in					
\cite[\S 2.1]{BaerStroh}
\begin{equation*} 
\begin{split}
&L^2(\Sigma_1, \mathscr{S}_{L,E}(\Sigma_1)) = \range{P_{\intervallro{a_1}{\infty}{}}(t_1)} \oplus 
\range{P_{\intervallo{-\infty}{a_1}{}}(t_1)}, \\
&L^2(\Sigma_2, \mathscr{S}_{L,E}(\Sigma_2)) = \range{P_{\intervallo{a_2}{\infty}{}}(t_2)} \oplus 
\range{P_{\intervalllo{-\infty}{a_2}{}}(t_2)}.
\end{split}
\end{equation*}
With respect to these decompositions, $Q(t_2,t_1), \tilde{Q}(t_2,t_1)$
are given by matrices 
\begin{equation*}
Q(t_2,t_1)=\left(\begin{matrix}
Q^{> a_2}_{\geq a_1}(t_2,t_1) & Q^{> a_2}_{< a_1}(t_2,t_1) \\[5pt]
Q^{\leq a_2}_{\geq a_1}(t_2,t_1) & Q^{\leq a_2}_{< a_1}(t_2,t_1)
\end{matrix}\right), \  
\tilde{Q}(t_2,t_1)=\left(\begin{matrix}
\tilde{Q}^{> a_2}_{\geq a_1}(t_2,t_1) & \tilde{Q}^{> a_2}_{< a_1}(t_2,t_1) \\[5pt]
\tilde{Q}^{\leq a_2}_{\geq a_1}(t_2,t_1) & \tilde{Q}^{\leq a_2}_{< a_1}(t_2,t_1)
\end{matrix}\right).
\end{equation*}
In particular, the individual entries of the block matrix $Q(t_2,t_1)$ are given by
\begin{equation}\label{qgenentries}
\begin{split}
Q^{> a_2}_{\geq a_1}(t_2,t_1) &:= P_{\intervallo{a_2}{\infty}{}}(t_2)\circ Q(t_2,t_1)\circ P_{\intervallro{a_1}{\infty}{}}(t_1), \\
Q^{\leq a_2}_{< a_1}(t_2,t_1) &:= P_{\intervalllo{-\infty}{a_2}{}}(t_2)\circ Q(t_2,t_1)\circ P_{\intervallo{-\infty}{a_1}{}}(t_1), \\
Q^{> a_2}_{< a_1}(t_2,t_1) &:= P_{\intervallo{a_2}{\infty}{}}(t_2)\circ Q(t_2,t_1)\circ P_{\intervallo{-\infty}{a_1}{}}(t_1), \\
Q^{\leq a_2}_{\geq a_1}(t_2,t_1) &:= P_{\intervalllo{-\infty}{a_2}{}}(t_2)\circ Q(t_2,t_1)\circ P_{\intervallro{a_1}{\infty}{}}(t_1), 
\end{split}
\end{equation}  
with similar identities for the matrix entries of $\tilde{Q}(t_2,t_1)$. 
For the special values $a_1=a_2=0$ we will also use the abbreviations
\begin{equation}\label{Qabbreviation}
\begin{array}{ccl}
Q_{++}(t_2,t_1) &=& Q^{> 0}_{\geq 0}(t_2,t_1), \\[5pt]
Q_{--}(t_2,t_1) &=& Q^{\leq 0}_{< 0}(t_2,t_1), \\
\end{array}
\quad 
\begin{array}{ccl}
Q_{+-}(t_2,t_1) &=&  Q^{> 0}_{< 0}(t_2,t_1), \\
Q_{-+}(t_2,t_1) &=&  Q^{\leq 0}_{\geq 0}(t_2,t_1).
\end{array}
\end{equation}
The crucial result by the first named author in \cite{OD2} is that
the diagonal entries $Q^{> a_2}_{\geq a_1}$ and $Q^{\leq a_2}_{< a_1}$ 
as well as $\tilde{Q}^{> a_2}_{\geq a_1}$ and $\tilde{Q}^{\leq a_2}_{< a_1}$ 
are $\Gamma$-Fredholm.
\begin{theo}[Theorem 9.2 in \cite{OD2}]\label{qgengammafred}
The following operators are $\Gamma$-Fredholm for any $a_1,a_2 \in \R$, namely
\begin{equation*}
\begin{split}
Q^{> a_2}_{\geq a_1}(t_2,t_1), \tilde{Q}^{> a_2}_{\geq a_1}(t_2,t_1) 
&\in \mathscr{F}_\Gamma(\range{P_{\geq a_1}(t_1)},\range{P_{> a_2}(t_2)}),\\
Q^{\leq a_2}_{< a_1}(t_2,t_1), \tilde{Q}^{\leq a_2}_{< a_1}(t_2,t_1) 
&\in \mathscr{F}_\Gamma(\range{P_{< a_1}(t_1)},\range{P_{\leq a_2}(t_2)}). 
\end{split}
\end{equation*}
with the various $\Gamma$-indices related by the following formulae
\begin{align*}
\Index_\Gamma\left(Q^{> a_2}_{\geq a_1}(t_2,t_1)\right) &= \Index_\Gamma\left(Q_{++}(t_2,t_1)\right) \\
&+\Chi_{\SET{a_2 > 0}}\dim_\Gamma\left(\range{P_{(0,a_2]}(t_2)}\right)-
\Chi_{\SET{a_2 < 0}}\dim_\Gamma\left(\range{P_{(a_2,0]}(t_2)}\right)\\
&+ \Chi_{\SET{a_1 < 0}}\dim_\Gamma\left(\range{P_{[a_1,0)}(t_1)}\right)-
\Chi_{\SET{a_1 > 0}}\dim_\Gamma\left(\range{P_{[0,a_1)}(t_1)}\right),\\
\Index_\Gamma\left(Q^{\leq a_2}_{< a_1}(t_2,t_1)\right) &= \Index_\Gamma\left(Q_{--}(t_2,t_1)\right)
\\ &+ \Chi_{\SET{a_2 < 0}}
\dim_\Gamma\left(\range{P_{(a_2,0]}(t_2)}\right)-\Chi_{\SET{a_2 > 0}}\dim_\Gamma\left(\range{P_{(0,a_2]}(t_2)}\right)\\
&+ \Chi_{\SET{a_1 > 0}}\dim_\Gamma\left(\range{P_{[0,a_1)}(t_1)}\right)-\Chi_{\SET{a_1 < 0}}
\dim_\Gamma\left(\range{P_{[a_1,0)}(t_1)}\right),
\end{align*}
where we abbreviated $\Chi_{\SET{a>0}} := \Chi_{(0,\infty)}(a)$ and $\Chi_{\SET{a<0}}:= \Chi_{(-\infty,0)}(a)$.
Exactly the same identities hold when $Q$ is replaced by $\tilde{Q}$.
\end{theo}

The $\Gamma$-Fredholmness of the Dirac operators $D^{E_L}_{\pm}$ under generalised (a)APS
boundary conditions has been shown by the first named author in \cite{OD2}.

\begin{theo}[Theorem 9.3 in \cite{OD2}]\label{indexDgaAPStwist}						
The operators \clef{gammafredoperres} 
are $\Gamma$-Fredholm with 
\begin{align*}
\Index_{\Gamma}(D^{E_L}_{+,\mathrm{APS}(a_1,a_2)})&= \Index_\Gamma\left(Q^{\leq a_2}_{< a_1}(t_2,t_1)\right), \
\Index_{\Gamma}(D^{E_L}_{+,\mathrm{aAPS}(a_1,a_2)})=\Index_\Gamma\left(Q^{> a_2}_{\geq a_1}(t_2,t_1)\right), \\
\Index_{\Gamma}(D^{E_L}_{-,\mathrm{APS}(a_1,a_2)})&= \Index_\Gamma\left(\tilde{Q}^{> a_2}_{\geq a_1}(t_2,t_1)\right), \
\Index_{\Gamma}(D^{E_L}_{-,\mathrm{aAPS}(a_1,a_2)})= \Index_\Gamma\left(\tilde{Q}^{\leq a_2}_{< a_1}(t_2,t_1)\right).
\end{align*}
\end{theo}

\section{Fundamentals of the Gamma spectral flow}\label{chap:Specf}

The spectral flow is used in \cite{BaerStroh} to connect the 
indices of the operators $Q^{\bullet}_{\bullet}(t_2,t_1)$ with geometric quantities. 
We introduce a definition of spectral flow in our setting and give an analytic 
expression at the end of this section. We refer to the general results of
\cite{bencarphiadfydwoj} for the case of a general von Neumann algebra 
and specify as well as modify their results to our case. 

\subsection{Recapitulation of the spectral flow}\label{chap:specflow-sec:gammaspecflow-1}

\subsubsection{Spectral flow following Phillips}

Let $I\subset \R$ be a closed interval and $\SET{S_t}_{t\in I}$ a continuous path of 
self-adjoint Fredholm operators with fixed domain $\mathscr{D}$ in a Hilbert space $\mathcal{H}$. 
The heuristical idea of the spectral flow is the net number of eigenvalues with multiplicities which 
pass through zero from negative to positive along the path. 
Let us begin with an analytic definition of the spectral flow for self-adjoint 
Fredholm operators, essentially due to Phillips \cite{phillips1996}.

\begin{defi}\label{defispectralflowordinary}
Let $\SET{S_t}_{t\in I}$ be a continuous path of 
self-adjoint Fredholm operators with domain $\mathscr{D}$ in a Hilbert space $\mathcal{H}$.
Choose an $N \in \N$ and a partition of the time interval $I=[t_1,t_2]$,
$$
t_1=\tau_0 < \tau_1 < ... < \tau_N=t_2,
$$ 
and positive real numbers 
$a_0,a_1,...,a_N$, such that $t\mapsto \Chi_{[-a_i,a_i]}(S_t)$ 
is a continuous family of finite-rank projection on $[\tau_{i-1},\tau_i]$ for each 
$i\in \SET{1,...,N}$. Then the spectral flow of the path is defined to be the number
\begin{equation}\label{defispectralflowordinaryformula}
\Specflow{S_t}{t\in I}:=\sum_{i=1}^N\Bigl[\dim\range{\Chi_{[0,a_i]}(S_{\tau_{i-1}})}-
\dim\range{\Chi_{[0,a_i]}(S_{\tau_{i}})}\Bigr].
\end{equation}
\end{defi}
The spectral flow is independent of the choice of the partition and 
the choice of real numbers outside of the spectrum. It is an integer-valued 
quantity, which adds up under concatenations of paths and is invariant 
under path homotopies with fixed endpoints in the space of 
self-adjoint Fredholm operators. See Theorem and Proposition 2 and 3 in \cite{phillips1996}
for bounded Fredholm operators. Extensions to densely defined self-adjoint Fredholm operators are 
intensively studied e.g. in \cite{leschspecflowunbound} and 
\cite{boossbavnbekleschphillips2005}. \medskip

\subsubsection{Spectral flow following Lesch}
We discuss another expression of the spectral flow 
formula. This notion uses relative indices, which can 
be generalized to the $\Gamma$-setting. We introduce the Calkin-map
\begin{equation*} 
\Uppi: \quad \mathscr{B}(\mathscr{H}) \quad \rightarrow \quad 
\quotspace{\mathscr{B}(\mathscr{H})}{\mathscr{K}(\mathscr{H})},
\end{equation*}
where $\mathscr{B}(\mathscr{H})$ denotes the space of bounded linear and
$\mathscr{K}(\mathscr{H})$ denotes the space compact operators on $\mathscr{H}$.
Denote with 
$$ 
\Index(P\vert Q):=\Index(P:\range{Q}\rightarrow \range{P}),
$$
the relative index of two projections $P$ and $Q$ on $\mathcal{H}$ such 
that $(P-Q)\in \mathscr{K}(\mathcal{H})$. Such a pair of projections is a 
\textit{Fredholm pair}; see e.g. \cite{avronseilersimon} for further conceptual details.  Abbreviate $P_{\geq \lambda}(t):= \Chi_{[\lambda, \infty)}(S_t)$. 
Relying on \cite[Sec.3]{leschspecflowunbound}, the existence of finite-rank projections in 
Definition \ref{defispectralflowordinary} can be relaxed to the condition that the path 
$\SET{P_{\geq \lambda}(t)}_{t\in I}$ of projections satisfies for all $t,t' \in [\tau_{j-1},\tau_j]$ 
and in each subinterval $[\tau_{j-1},\tau_j]$ of the partition
\begin{equation}
\norm{\Uppi \Bigl( P_{\geq a_j}(t) \Bigr) - \Uppi \Bigl( P_{\geq a_j}(t') \Bigr)}{\mathscr{B}(\mathcal{H})}< 1,
\end{equation}
where, by a slight abuse of notation, the norm $\norm{\cdot}{\mathscr{B}(\mathscr{H})}$ denotes the 
canonical induced norm on the quotient 
$\quotspace{\mathscr{B}(\mathscr{H})}{\mathscr{K}(\mathscr{H})}$.
In particular, \cite[Thm. 3.6]{leschspecflowunbound} and \cite[Cor. 3.8]{leschspecflowunbound} 
motivates the following definition of spectral flow. 

\begin{defi}\label{defispectralflowordalt}
Let $\SET{S_t}_{t\in I}$ be a continuous path of 
self-adjoint Fredholm operators with domain $\mathscr{D}$ in a Hilbert space $\mathcal{H}$. 
Choose a subdivision $t_1=\tau_0 < \tau_1 < ... < \tau_N=t_2$ of $I$ fine enough such that 
for $t,t'\in [\tau_{j-1},\tau_j]$ and any $j$
$$
\norm{\Uppi \Bigl( P_{\geq 0}(t) \Bigr) - \Uppi \Bigl( P_{\geq 0}(t') \Bigr)}{\mathscr{B}(\mathcal{H})}< 1,
$$
where, as before, by a slight abuse of notation 
the norm $\norm{\cdot}{\mathscr{B}(\mathscr{H})}$ denotes the 
canonical induced norm on the Calkin algebra
$\quotspace{\mathscr{B}(\mathscr{H})}{\mathscr{K}(\mathscr{H})}$.
Then the spectral flow of the path is defined to be the number
\footnote{In \cite{leschspecflowunbound} the order of projections is reversed due to a reversed mapping orientation in the definition of the relative index,  leading
to a different sign convention.}
\begin{equation}\label{defispectralflowordinaryformulaalt}
\Specflow{S_t}{t\in I}:=\sum_{i=1}^N\Index\bigl(P_{\geq 0}(\tau_{j-1})\,\vert\,
P_{\geq 0}(\tau_{j})\bigr).
\end{equation}
\end{defi}
This definition is going to be the key ingredient in defining a spectral flow in 
von Neumann algebras and in particular for our $\Gamma$-setting. If there 
exist positive real numbers $a_0,a_1,....,a_N$ which are neither in the spectrum 
nor their intermediate intervals $[a_{j-1},a_j]$ have non-trivial intersections with 
the essential spectrum for $t$ in any subinterval $[t_{j-1},t_j]$, formula 
\clef{defispectralflowordinaryformula} is recovered, 
as the family of projectors $\SET{P_{>a_j}}$ becomes continuous and 
$$ \Index\bigl(P_{\geq 0}(\tau_{j-1})\,\vert\,P_{\geq 0}(\tau_{j})\bigr)=
\dim\range{P_{[0,a_i]}(\tau_{j-1})}-\dim\range{P_{[0,a_i]}(\tau_{j})}.$$

\subsubsection{Spectral flow and the eta-invariant}

Let $A$ be an elliptic and self-adjoint (pseudo-)differential operator of 
order $d\geq 0$ on a closed manifold $\Sigma$ of dimension $n$. In particular, $A$ admits 
discrete spectrum $\upsigma(A)=\SET{\lambda_j\,\vert\, j\in \N_0}$ and we may define the
associated \textit{eta-invariant} as in \cite{APSIII}
\begin{equation*}
\eta(z;A):=\sum_{\lambda \in \upsigma(A)\setminus \SET{0}} \frac{\mathrm{sign}(\lambda)}{\vert \lambda \vert^z}
= \frac{1}{\Gamma(\frac{z+1}{2})}\int_0^\infty s^{\frac{z-1}{2}} \Tr{}{A\expe{-sA^2}}\differ s,
\quad \Rep{z} > \frac{n}{d}.
\end{equation*}

It is holomorphic for $\Rep{z}>\frac{n}{d}$ and extends meromorphically to $\C$ 
with at most simple poles at $z\in (-1/d)\Z\setminus\SET{0}$.
The residue of $\eta(z;A)$ at $z=0$ is determined by the complete symbol of the operator $A$ (see \cite[Prop.2.8]{APSIII}), implying that it is in general non-zero. However, for $n$ odd, the residue vanishes (see \cite[Thm.4.5]{APSIII}); same question has been studied on even-dimensional manifolds in \cite{GILKEY1981290} (Theorem 0.1, Lemma 1.1/2). In those cases, $z=0$ is a removeable singularity and $\eta(0;A)$ the well-defined \textit{eta-invariant}
\begin{equation}\label{etainv}
\eta(A):=\eta(0,A)=\frac{1}{\sqrt{\uppi}} \int_0^\infty s^{-\frac{1}{2}} \Tr{}{A\expe{-sA^2}}\differ s.
\end{equation}
The \textit{reduced eta-invariant} (see \cite{APSII}), also called the \textit{xi-invariant}, 
arises in the APS index Theorem \cite[Thm.3.10]{APSI} for compact manifolds with boundary and is defined by
\begin{equation}\label{xiinv}
\xi(A):=\frac{\eta(A)+\dim\kernel{A}}{2}.
\end{equation}
Its connection to the spectral flow 
is studied in \cite{APSIII}. The following fact sums up their result.
\begin{prop}[cf.  Prop. 1.1 in \cite{daixhang}]\label{propspecflow}
Let $\SET{S_t}_{t\in[t_1,t_2]}$ be a smooth family of elliptic, self-adjoint operators of fixed order $d$ on a 
closed Riemannian manifold $\Sigma$ of odd dimension. Then the family of the associated 
eta-invariants $t\mapsto \eta(S_t)$ is piecewise smooth and admits a decomposition of the form
\begin{equation}\label{specflowformulaI}
\xi(S_{t_2})-\xi(S_{t_1})=\Specflow{S_t}{t\in[t_1,t_2]}+\int_{t_1}^{t_2}\frac{\differ}{\differ t}\xi(S_t)\differ t
\end{equation}
where the spectral flow is piecewise constant and the integrand
is smooth in $t$.
\end{prop} 

We should emphasize, that while $\frac{\differ}{\differ t}\xi(S_t)$ is smooth in $t$, 
$\xi(S_t)$ is not. In fact $\xi(S_t)$ is only piecewise smooth in $t$, and hence the integral
$\int_{t_1}^{t_2}\frac{\differ}{\differ t}\xi(S_t)\differ t$ does not equal $\xi(S_{t_2}) - \xi(S_{t_1})$.

\subsection{$\Gamma$-Fredholm pairs and their $\Gamma$-indices}\label{chap:specflow-sec:gammafredpair}

We transfer concepts of relative Fredholm theory 
and relative index to the $\Gamma$-setting. 
The discussion is based on \cite{bencarphiadfydwoj} and \cite{carphil} for general semifinite 
von Neumann algebras $\mathscr{N}$. Our special case here is 
$\mathscr{N}=\mathscr{B}_\Gamma(\mathscr{H})$, the space of bounded linear 
operators on a free or projective Hilbert $\Gamma$-module $\mathscr{H}$.
We introduce the $\Gamma$-Calkin-map
\begin{equation*} 
\Uppi_\Gamma: \quad \mathscr{B}_\Gamma(\mathscr{H}) \quad \rightarrow \quad \quotspace{\mathscr{B}_\Gamma(\mathscr{H})}{\mathscr{K}_\Gamma(\mathscr{H})},
\end{equation*}
where $\mathscr{K}_\Gamma(\mathscr{H})$ denotes the space of $\Gamma$-compact operators,
defined e.g. in \cite[Definition 3.9]{OD2}. An operator $A \in \mathscr{B}_\Gamma(\mathscr{H})$ is 
$\Gamma$-Fredholm, we write $A \in \mathscr{F}_\Gamma(\mathscr{H})$, if and only if its image
under the $\Gamma$-Calkin map is invertible.  The definition holds for 
operators between different Hilbert $\Gamma$-modules with a parallel notation. 

\begin{defi}\label{relative-F-pair-defi}
\begin{enumerate}
\item Given two projectors $P$ and $Q$ in (free, projective) Hilbert $\Gamma$-modules $\mathscr{H}_1$ $\mathscr{H}_2$ such that $\range{Q}\subset \mathscr{H}_1$ and $\range{P}\subset\mathscr{H}_2$, an operator $A \in \mathscr{B}_\Gamma(\range{Q},\range{P})$ is called $(P,Q)\text{-}\Gamma$-Fredholm if $A \in \mathscr{F}_\Gamma(\range{Q},\range{P})$. The corresponding $\Gamma$-index is
\begin{equation*}
\Index_{\Gamma}\left(A\vert_{\range{Q}\rightarrow\range{P}}\right)
= \dim_\Gamma(\kernel{A}\cap Q(\mathscr{H}_1))-\dim_\Gamma(\kernel{A^\ast}\cap P(\mathscr{H}_2)).
\end{equation*}
\item If $\mathscr{H}_1 = \mathscr{H_2}$ and 
$PQ$ is $(P,Q)\text{-}\Gamma$-Fredholm, i.e. $PQ \in \mathscr{F}_\Gamma(\range{Q},\range{P})$,
we call $(P,Q)$ a $\Gamma$\textit{-Fredholm pair}, with the corresponding $\Gamma$-index
abbreviated as
\begin{equation}\label{upgammaindexfredpair}
\Index_\Gamma(P\vert Q):= \Index_\Gamma(P\vert_{\range{Q}\rightarrow \range{P}})
\equiv \Index_\Gamma(PQ\vert_{\range{Q}\rightarrow \range{P}}).
\end{equation}
\end{enumerate}
\end{defi}

An example of such operators are the restricted projections with $P=P_{\blacktriangledown b}$ and $Q=P_{\blacktriangle a}$ and Hilbert $\Gamma$-modules $\mathscr{H}_1=\mathscr{H}_2=L^2$ from \cite[Lemma 7.1]{OD2}. The ordinary properties of $\Gamma$-Fredholm operators and their $\Gamma$-indices carry over to this specific concept. We refer to \cite[Sec.4]{bencarphiadfydwoj} for details. A condition for $(P,Q)$ to form
a $\Gamma$-Fredholm pair is stated in \cite[Lemma 4.1]{bencarphiadfydwoj} as follows.

\begin{lem}\label{gammacalkfredpairproj}
Projections $P,Q \in \mathscr{B}_\Gamma(\mathscr{H})$
form a $\Gamma$-Fredholm pair if and only if
\begin{equation*}
\norm{\Uppi_\Gamma(P)-\Uppi_\Gamma(Q)}{\mathscr{B}(\mathscr{H})} < 1,
\end{equation*}
where, by a slight abuse of notation, the norm $\norm{\cdot}{\mathscr{B}(\mathscr{H})}$ denotes the 
canonical induced norm on the quotient 
$\quotspace{\mathscr{B}_\Gamma(\mathscr{H})}{\mathscr{K}_\Gamma(\mathscr{H})}$.
\end{lem}

The $\Gamma$-index of a $\Gamma$-Fredholm pair has the following properties.

\begin{prop}\label{propindexfredpair} %
For orthogonal projections $P,Q,R \in \mathscr{B}_\Gamma(\mathscr{H})$ the following holds:
\begin{itemize}
\item[(1)]
if $PQ \in \mathscr{F}_\Gamma(\range{Q},\range{P})$, then $\dim_\Gamma\kernel{P-Q\pm \Iop{}} < \infty$ and 
\begin{equation}\label{upgammaindexfredpairalt}
\Index_\Gamma(P\vert Q)= \dim_\Gamma\kernel{P-Q+ \Iop{}}-\dim_\Gamma\kernel{P-Q- \Iop{}}.
\end{equation}
\item[(2)] $PQ \in \mathscr{F}_\Gamma(\range{Q},\range{P})$ if and only if 
$QP \in \mathscr{F}_\Gamma(\range{P},\range{Q})$ with 
$$
\Index_\Gamma(Q\vert P)=-\Index_\Gamma(P\vert Q).
$$
\item[(3)] The duals $P^\perp = \Iop{} - P$ and $Q^\perp = \Iop{} - Q$ satisfy 
$P^\perp Q^\perp \in \mathscr{F}_\Gamma(\range{Q^\perp},\range{P^\perp})$ if and only if 
$PQ \in \mathscr{F}_\Gamma(\range{Q},\range{P})$ with 
$$
\Index_\Gamma(P^\perp\vert Q^\perp)=-\Index_\Gamma(P\vert Q).
$$
\item[(4)] If $PR \in \mathscr{F}_\Gamma(\range{R},\range{P})$ and 
$RQ \in \mathscr{F}_\Gamma(\range{Q},\range{R})$ with
\begin{equation*}
\norm{\Uppi_\Gamma(R)-\Uppi_\Gamma(Q)}{\mathscr{B}(\mathscr{H})} < \frac{1}{2} \quad \text{and} \quad \norm{\Uppi_\Gamma(P)-\Uppi_\Gamma(R)}{\mathscr{B}(\mathscr{H})} < \frac{1}{2} \quad,
\end{equation*}
then $(PQ) \in \mathscr{F}_\Gamma(\range{Q},\range{P})$ with 
$$
\Index_\Gamma(P\vert Q)=\Index_\Gamma(P\vert R)+\Index_\Gamma(R\vert Q).
$$
\item[(5)] If $\norm{P-Q}{\mathscr{B}(\mathscr{H})} < 1$, then $\Index_\Gamma(P\vert Q)=0$.
\item[(6)] If $PQ \in \mathscr{F}_\Gamma(\range{Q},\range{P})$ and $U:\mathscr{H}\rightarrow \mathscr{H}$ is a unitary $\Gamma$-isomorphism, then $(UPU^{-1})\circ(UQU^{-1})$ is $\Gamma$-Fredholm pair with
$$
\Index_\Gamma(UPU^{-1}\vert UQU^{-1})=\Index_\Gamma(P\vert Q).
$$
\end{itemize}
\end{prop}

\begin{proof} 
Statement (1) can be proven purely algebraically as in \cite[Prop.3.1]{avronseilersimon} by rewriting the kernels,
such that \clef{upgammaindexfredpairalt} follows directly from \clef{upgammaindexfredpair}.
Indeed,  note that for all $x\in\range{Q}$ and all $y\in\range{P}$
\begin{eqnarray*}
Px=0 \quad \Leftrightarrow \quad &Px= Qx-x \quad \Leftrightarrow \quad &(P-Q+\Iop{})x=0, \\
Qy=0 \quad \Leftrightarrow \quad &Qy= Py-y \quad \Leftrightarrow \quad &(P-Q-\Iop{})y=0.
\end{eqnarray*}
From here we conclude
\begin{eqnarray*}
\kernel{P\vert_{\range{Q}\rightarrow\range{P}}} &=& 
\SET{x\in\range{Q}\,\vert\, Px=0}=
\kernel{P-Q+\Iop{}},\\
\kernel{Q\vert_{\range{P}\rightarrow\range{Q}}} &=& 
\SET{y\in\range{P}\,\vert\, Qy=0}=
\kernel{P-Q-\Iop{}}. 
\end{eqnarray*}
Since $PQ \in \mathscr{F}_\Gamma(\range{Q},\range{P})$, the range of 
$P\vert_{\range{Q}\rightarrow\range{P}}$ is closed and 
\begin{eqnarray*}
\dim_\Gamma\kernel{P-Q+\Iop{}} &=& \dim_\Gamma\kernel{P\vert_{\range{Q}\rightarrow\range{P}}} < \infty \\
\dim_\Gamma\kernel{P-Q-\Iop{}} &=& \dim_\Gamma\kernel{Q\vert_{\range{P}\rightarrow\range{Q}}},
\end{eqnarray*}
such that the index formula \clef{upgammaindexfredpairalt} follows from Definition \ref{relative-F-pair-defi}.

As the condition for Fredholmness in Lemma \ref{gammacalkfredpairproj} is symmetric in $P$ and $Q$, 
the equivalence (2) follows easily. The statement (3) follows from the fact that 
the linear $\Gamma$-Calkin map gives $\Uppi_\Gamma(\Iop{})=\Iop{}$ and hence
\begin{equation*}
\norm{\Uppi_\Gamma(P^\perp)-\Uppi_\Gamma(Q^\perp)}{\mathscr{B}(\mathscr{H})}=\norm{\Iop{}-\Uppi_\Gamma(P)-\Iop{}+\Uppi_\Gamma(Q)}{\mathscr{B}(\mathscr{H})}=\norm{\Uppi_\Gamma(P)-\Uppi_\Gamma(Q)}{\mathscr{B}(\mathscr{H})}.
\end{equation*}
The $\Gamma$-index relations in (2) and (3) follow from \clef{upgammaindexfredpairalt}.
The statements (4) and (5) are proven as in \cite[Lemma 4.3 and Remark 4.4]{bencarphiadfydwoj} 
for general semifinite von Neumann algebras. 
Statement (6) is also a consequence of (1): the kernels in \clef{upgammaindexfredpairalt} can be rewritten as follows:
\begin{equation*}
\kernel{UPU^{-1}-UQU^{-1}\pm \Iop{}}=\kernel{U(P-Q\pm \Iop{})U^{-1}}=U(\kernel{P-Q\pm \Iop{}}).
\end{equation*}
Because $U$ is a unitary $\Gamma$-isomorphism, we obtain
$$
\kernel{UPU^{-1}-UQU^{-1}\pm \Iop{}}\cong\kernel{P-Q\pm \Iop{}}
$$
and the $\Gamma$-dimensions as well as the relative $\Gamma$-indices coincide. 
\end{proof}

\subsection{Algebraic definition of the $\Gamma$-spectral flow}
\label{chap:specflow-sec:gammaspecflow-2}

Our definition is an analogue of Definition \ref{defispectralflowordalt}.
It mimics the approach for spectral flow in any semifinite von Neumann algebra by
\cite{carphil} and \cite{bencarphiadfydwoj},  adapted to our setting.

\begin{defi}
Let $\SET{S_t}_{t\in I}$ be a continuous path of 
self-adjoint $\Gamma$-Fredholm operators with domain $\mathscr{D}$ in a Hilbert space $\mathcal{H}$,
both Hilbert $\Gamma$-modules. Choose a sufficiently fine subdivision of $I=[t_1,t_2]$, 
$$
t_1=\tau_0 < \tau_1 < ... < \tau_N=t_2,
$$ 
such that 
for $t,t'\in [\tau_{j-1},\tau_j]$ and any $j$ the corresponding spectral projections 
$P_{\geq \lambda}(t):= \Chi_{[\lambda, \infty)}(S_t) \in \mathscr{B}_\Gamma(\mathscr{H})$
satisfy
\begin{equation*}
\norm{\Uppi_\Gamma(P_{\geq 0}(t))-\Uppi_\Gamma(P_{\geq 0}(t'))}{\mathscr{B}(\mathscr{H})} < 1.
\end{equation*}
The $\Gamma$-spectral flow of the path $\SET{S_t}_{t\in I}$ is then defined by
\begin{equation}\label{gammaspectralflow}
\specflow{S_t}{t\in I}:=\sum_{j=1}^N\Index_\Gamma
\Bigl(P_{\geq 0}(\tau_{j-1})\vert P_{\geq 0}(\tau_{j})\Bigr).
\end{equation}
\end{defi}

\begin{rem}
\begin{itemize}
\item[(i)] The definition is independent of the chosen partition and works 
equally well for all von Neumann algebras of type I, II and mixed; see \cite[Thm.2.1]{bencarphiadfydwoj}. 
The $\Gamma$-spectral flow is a priori a real non-integer number due to the fact that the 
$\Gamma$-indices are in 
general real-valued. \medskip

\item[(ii)] If the path of projections $\SET{P_{\geq 0}(t)}$ is continuous, we can choose a fine 
enough partition such that for all $t,s\in[\tau_{j-1},\tau_j]$ and any $j$
\begin{equation*}
\norm{P_{\geq 0}(t))-P_{\geq 0}(s))}{\mathscr{B}(\mathscr{H})} < 1.
\end{equation*}
The statement in Proposition \ref{propindexfredpair} (5) then implies that 
$\Index_\Gamma(P_{\geq 0}(r)\vert P_{\geq 0}(s))=0$ 
for each $r,s\in[\tau_{j-1},\tau_j]$ in the partition. As a consequence $\specflow{S_t}{t\in I}$ 
vanishes. This fits with the heuristical expectation that the spectral flow is only non-trivial 
if the path of projections has discontinuities. \medskip

\item[(iii)] If $\SET{S_t}_{t\in[t_1,t_2]}$ and $\SET{\bar{S}_t}_{t\in[t_1,t_2]}$ are two 
homotopic continuous paths with $\bar{S}_{t_1}=S_{t_1}$ and $\bar{S}_{t_2}=S_{t_2}$, 
the $\Gamma$-spectral flows of these paths coincide. Hence the $\Gamma$-spectral flow is a 
homotopy invariant. This is proven for general von Neumann algebras of type I and II in 
\cite[Prop.4.7]{bencarphiadfydwoj}.  
\end{itemize}
\end{rem}

\subsection{$\Gamma$-spectral flow with variable Hilbert spaces}

So far we considered the case that the domain $\mathscr{D}$ and 
Hilbert space $\mathscr{H}$ of each $S_t$ are the same for each $t \in [t_1,t_2]$. 
Below, we apply the $\Gamma$-spectral flow concept to the path of 
$\Gamma$-invariant hypersurface Dirac operators with \emph{variable} 
Hilbert space structures. Assume, $\SET{\mathscr{H}(t)}_{t\in [t_1,t_2]}$ 
is a family of Hilbert $\Gamma$-modules with unitary $\Gamma$-isomorphisms
\begin{equation}\label{isometrybundlehilbertgammamodules}
Q(\tau_2,\tau_1)\,:\,\mathscr{H}(\tau_1)\rightarrow\mathscr{H}(\tau_2),
\end{equation}
for any $\tau_1, \tau_2 \in [t_1,t_2]$, continuous in both time entries.
The notation corresponds to \eqref{Q-def}, since these will be precisely those
unitary transformations that we will use for $S_t := A_t^{E_L}$.
Let $\SET{S_t}_{t\in [t_1,t_2]}$ be a family of 
self-adjoint $\Gamma$-Fredholm operators with each $S_\tau$
having dense domain $\mathscr{D}(\tau) \subset \mathscr{H}(\tau)$,
such that 
$$
Q(\tau_2,\tau_1) \mathscr{D}(\tau_1) = \mathscr{D}(\tau_2).
$$
As in \cite{ronge} or \cite{rongedendungen} we introduce the \textit{evolved operator}
\begin{equation*}
\hat{S}_t:=Q(t_1,t) S_t\,Q(t,t_1), \quad \mathscr{D}(\hat{S}_t) = \mathscr{D}(t_1).
\end{equation*}   
By construction, $\SET{\hat{S}_t}_{t\in [t_1,t_2]}$ is a family of 
self-adjoint $\Gamma$-Fredholm operators with fixed domain $\mathscr{D}(t_1)$ in a 
fixed Hilbert space $\mathscr{H}(t_1)$,  with spectral projections 
$P_{J}(t) = \Chi_{J}(S_t)$ and $\hat{P}_{J}(t) = \Chi_{J}(\hat{S}_t)$
related for any $J \subset \R$ by
\begin{equation}\label{PPhat}
\begin{split}
&\hat{P}_J(t)=Q(t_1,t)P_J(t)Q(t,t_1), \\
&\range{\hat{P}_J(t)}=Q(t_1,t)(\range{P_J(t)}).
\end{split}
\end{equation}
Assuming that $\SET{\hat{S}_t}_{t\in [t_1,t_2]}$ is a continuous path,
we may consider the $\Gamma$-spectral flow with 
respect to this path.  This leads to the following definition.

\begin{defi}\label{defimodgammaspectralflow}
Consider a family $\SET{S_t : \mathscr{D}(t) \rightarrow \mathscr{H}(t)}_{t\in I}$ of self-adjoint $\Gamma$-Fredholm operators and the corresponding family $\SET{\hat{S}_t:\mathscr{D}(t_1)\rightarrow \mathscr{H}(t_1)}_{t\in I}$ of evolved operators.  Assume that the evolved operators
$\SET{\hat{S}_t}_{t\in I}$ form a continuous path in
$\mathscr{F}_\Gamma^{\text{sa}}(\mathscr{D}(t_1),\mathscr{H}(t_1))$.
Choose a sufficiently fine subdivision of $I=[t_1,t_2]$
$$
t_1=\tau_0 < \tau_1 < ... < \tau_N=t_2,
$$ 
such that 
for $t,t'\in [\tau_{j-1},\tau_j]$ and any $j$ the corresponding spectral projections 
$\hat{P}_{\geq \lambda}(t):= \hat{P}_{[\lambda,\infty)}(t) \in 
\mathscr{B}_\Gamma(\mathscr{H}(t_1))$
satisfy
\begin{equation*}
\norm{\Uppi_\Gamma(\hat{P}_{\geq 0}(t))-
\Uppi_\Gamma(\hat{P}_{\geq 0}(t'))}{\mathscr{B}(\mathscr{H}(t_1))} < 1.
\end{equation*}
The $\Gamma$-spectral flow of the path $\SET{S_t}_{t\in I}$ is then defined by
\begin{equation}\label{modgammaspectralflow}
\specflow{S_t}{t\in I}:=\specflow{\hat{S}_t}{t\in I} \equiv \sum_{j=1}^N\Index_\Gamma
\Bigl(\hat{P}_{\geq 0}(\tau_{j-1})\Big\vert \hat{P}_{\geq 0}(\tau_{j})\Bigr).
\end{equation}
\end{defi}

It remains to prove that $\specflow{S_t}{t\in I}$ is well-defined,  i.e.  that the definition does not depend on the choice of the evolution operators
$Q(\tau_2,\tau_1)$ in \eqref{isometrybundlehilbertgammamodules},  which define the evolved operators and projections.

\begin{prop}\label{sf-sf-modified-formal} 
The spectral flow $\specflow{S_t}{t\in I}$ in Definition \ref{defimodgammaspectralflow} 
is well-defined,  i.e.  does not depend on the choice of the evolution operators
$Q(\bullet, \bullet)$ in \clef{isometrybundlehilbertgammamodules}.
\end{prop}
\begin{proof}
Suppose that there is another unitary $\Gamma$-isomorphism $U(\tau_2,\tau_1)\,:\,\mathscr{H}(\tau_1)\rightarrow\mathscr{H}(\tau_2)$ for any $\tau_1, \tau_2 \in [t_1,t_2]$,  next to $Q$,  
having the same properties as \clef{isometrybundlehilbertgammamodules}. We define a family of evolved operators with respect to $U$
$$
\check{S}_t:=U(t_1,t)S_tU(t,t_1),
$$ 
which is assumed to be a continous familiy of self-adjoint $\Gamma$-Fredholm operators in 
$\mathscr{H}(t_1)$ as well.  Setting $\check{P}_{\geq 0}(t):=U(t_1,t)P_{\geq 0}(t)U(t,t_1)$,  
the claim which we now need to prove is equality of the two following spectral flows
\begin{equation}\label{sf12}
\begin{split}
&\specflow{\hat{S}_t}{t\in I}=\sum_{j=1}^N\Index_\Gamma\Bigl(\hat{P}_{\geq 0}(\tau_{j-1})\Big\vert \hat{P}_{\geq 0}(\tau_{j})\Bigr),  \\
&\specflow{\check{S}_t}{t\in I} =\sum_{j=1}^N \Index_\Gamma
\Bigl(\check{P}_{\geq 0}(\tau_{j-1})\Big\vert \check{P}_{\geq 0}(\tau_{j})\Bigr).
\end{split}\end{equation}
We define $R(t_1):=Q(t_1,t)U(t,t_1)$,  which is a unitary $\Gamma$-isomorphism from $\mathscr{H}(t_1)$ to $\mathscr{H}(t_1)$ with $R^\ast(t_1)=U(t_1,t)Q(t,t_1)$. 
Then $\check{P}_{\geq 0}(t)$
is related to the evolved projector $\hat{P}_{\geq 0}(t)$ for each $t$ by
\begin{eqnarray*}
\check{P}_{\geq 0}(t)&=& U(t_1,t)P_{\geq 0}(t)U(t,t_1)\\
&=&U(t_1,t)Q(t,t_1)Q(t_1,t)P_{\geq 0}(t)Q(t,t_1)Q(t_1,t)U(t,t_1)\\
&=&\left(U(t_1,t)Q(t,t_1)\right)\circ \hat{P}_{\geq 0}(t)\circ\left(Q(t_1,t)U(t,t_1)\right) \\
&=&R^\ast(t_1)\hat{P}_{\geq 0}(t)R(t_1).
\end{eqnarray*}
Due to this unitary conjugation it is clear that a $\Gamma$-Fredholm pair of evolved projections with respect to $Q$ is also a $\Gamma$-Fredholm pair of evolved projections with respect to $U$.  \Cref{propindexfredpair} (6) then implies 
$$
\Index_\Gamma\Bigl(\hat{P}_{\geq 0}(\tau_{j-1})
\Big\vert \hat{P}_{\geq 0}(\tau_{j})\Bigr)
= \Index_\Gamma
\Bigl(\check{P}_{\geq 0}(\tau_{j-1})\Big\vert \check{P}_{\geq 0}(\tau_{j})\Bigr).
$$
This implies equality of the two sums in \eqref{sf12} and hence the claim.
\end{proof}
We will use this result later with $Q(t,t_1)$ in \eqref{Q-def},  and $U(t,t_1)$ being the natural isometry, induced by parallel transport, which is also a unitary $\Gamma$-isomorphism between $L^2$-section for different slices $\Sigma_t$ and in addition differentiable in time arguments.

\section{Gamma spectral flow and the $L^2$-Gamma eta-invariant}\label{chap:specflow-sec:gammaspecflow-3}

M. Ramachandran introduced in \cite{ramachandran} the $\Upgamma$-analogue 
for the eta function, replacing the trace with the $\Upgamma$-trace
for any  $\Upgamma$-invariant and geometric Dirac operator $A$. More precisely, we define
for any $\epsilon > 0$ the so-called truncated $\Upgamma$-\textit{eta functions} 
\begin{equation}\label{gammaetafun}
\begin{split}
&\eta^{<\epsilon}_\Upgamma(z;A) := \frac{1}{\Gamma(\frac{z+1}{2})}
\int_0^\epsilon s^{\frac{z-1}{2}} \Tr{\Upgamma}{A\expe{-sA^2}}\differ s, \\
&\eta^{>\epsilon}_\Upgamma(z;A) := \frac{1}{\Gamma(\frac{z+1}{2})}
\int_\epsilon^\infty s^{\frac{z-1}{2}} \Tr{\Upgamma}{A\expe{-sA^2}}\differ s.
\end{split}
\end{equation}
The former is well-defined for $\Re(z)> (n+1)/2$ due to the short time 
asymptotics of the integrand and is holomorphic at $z=0$ due to the
\textit{Bismut-Freed-cancellation property} 
(see \cite[Lem.3.1.1]{ramachandran} or \cite{bismutfreed})
\begin{equation}\label{bismutfreedgamma}
s^{-1/2} \ \Tr{\Upgamma}{A\expe{-sA^2}}=\mathcal{O}(1), \quad s\to 0.
\end{equation}
The latter function, $\eta^{>\epsilon}_\Upgamma(z;A)$ is bounded near $z=0$
due to the estimate in \cite[(3.1.17)]{ramachandran}. We can therefore evaluate 
both functions at $z=0$ and write
$$
\eta^{<\epsilon}_\Upgamma(A):= \eta^{<\epsilon}_\Upgamma(0,A), \quad 
\eta^{>\epsilon}_\Upgamma(A):= \eta^{>\epsilon}_\Upgamma(0,A),
$$
where $\eta^{<\epsilon}_\Upgamma(A)$,
$\eta^{>\epsilon}_\Upgamma(A)$ are referred to as 
the \textit{upper and lower truncated $\Upgamma$-eta invariants}, respectively.
We can now define for any $\epsilon > 0$
 the $\Upgamma$-\textit{eta invariant} as
\begin{equation}\label{gammaetainv}
\begin{split}
\eta_\Upgamma(A) &:= \frac{1}{\sqrt{\uppi}} 
\int_0^{\epsilon} s^{-\frac{1}{2}} \Tr{\Upgamma}{A\expe{-sA^2}}\differ s
+ \frac{1}{\sqrt{\uppi}}
\int_{\epsilon}^\infty s^{-\frac{1}{2}} \Tr{\Upgamma}{A\expe{-sA^2}}\differ s
\\[3mm] &\equiv \eta^{<\epsilon}_\Upgamma(A) + \eta^{>\epsilon}_\Upgamma(A).
\end{split}
\end{equation}
We also define the \textit{upper and lower truncated $\Upgamma$-xi invariants}
$$
\xi^{<\epsilon}_\Upgamma(A):=\frac{\eta^{<\epsilon}_\Upgamma(A)+\dim_\Upgamma\kernel{A}}{2},
\quad \xi^{>\epsilon}_\Upgamma(A):=\frac{\eta^{>\epsilon}_\Upgamma(A)+\dim_\Upgamma\kernel{A}}{2}. 
$$
We also define the $\Upgamma$-\textit{xi invariant}
\begin{equation}\label{gammaxi}
\xi_\Upgamma(A):=\frac{\eta_\Upgamma(A)+\dim_\Upgamma\kernel{A}}{2}
= \xi^{<\epsilon}_\Upgamma(A) + \xi^{>\epsilon}_\Upgamma(A).
\end{equation}

\subsection{Analytic expression for the $\Gamma$-spectral flow}
\label{chap:specflow-sec:gammaspecflow-4}

We want an analytic expression of the $\Gamma$-spectral flow 
$\specflow{A_t}{t\in I}$ of the continuous family of twisted hypersurface Dirac operators 
$A_t := A^{E_L}_t$ by means of the $\Gamma$-eta invariant as in Proposition \ref{propspecflow}.  Consider 
$$
U(t,t_1): L^2(\Sigma_{t_1}, \mathscr{S}_{L,E}(\Sigma_{t_1}))
\to L^2(\Sigma_t, \mathscr{S}_{L,E}(\Sigma_t)),
$$ 
the natural isometry of Hilbert $\Upgamma$-modules,  defined by parallel transport.  
In this subsection,  we will specifically use $U(t,t_1)$ instead of $Q(t,t_1)$ from \eqref{Q-def},  since the latter is not necessarily differentiable in $t$. 
The evolved operators are then defined by
\begin{equation}
\check{A}_t:=U(t_1,t)A_t\,U(t,t_1), \quad 
\expe{-s\check{A}^2_t}=U(t_1,t)\expe{-sA^2_t}U(t,t_1).
\end{equation}
Based on the same algebraic definition of spectral flow, as in \eqref{gammaspectralflow}
an analytic spectral flow formula has been proven in \cite[Prop.3.3]{AzWahl} for any type II von Neumann algebra setting.  Their result reads for any $\epsilon > 0$ as follows
\begin{equation}\label{gammaspectrunc2}
\begin{split}
\specflow{\check{A}_t}{t \in [t_1,t_2]} &= \xi^{>\epsilon}_\Gamma(\check{A}_{t_2})-\xi^{>\epsilon}_\Gamma(\check{A}_{t_1})+\sqrt{\frac{\epsilon}{\uppi}}\int_{t_1}^{t_2}\Tr{\Gamma}{\dot{\check{A}}_t\expe{-\epsilon \check{A}^2_t}} \differ t \\
&= \xi_\Gamma(\check{A}_{t_2})-\xi_\Gamma(\check{A}_{t_1})+\frac{\eta^{<\epsilon}_\Gamma(\check{A}_{t_1})-
\eta^{<\epsilon}_\Gamma(\check{A}_{t_2})}{2} \\ &+\sqrt{\frac{\epsilon}{\uppi}}\int_{t_1}^{t_2}
\Tr{\Gamma}{\dot{\check{A}}_t\expe{-\epsilon \check{A}^2_t}} \differ t.
\end{split}
\end{equation}  
We first observe,  that in this equality all instances of
$\check{A}_t$ can be replaced by $A_t$.  Namely,  we have the following

\begin{prop}\label{gammaspectrunc1-prop}
\begin{equation}\label{gammaspectrunc1}
\begin{split}
\specflow{A_t}{t \in [t_1,t_2]} &= \xi_\Gamma(A_{t_2})-\xi_\Gamma(A_{t_1})+\frac{\eta^{<\epsilon}_\Gamma(A_{t_1})-
\eta^{<\epsilon}_\Gamma(A_{t_2})}{2} \\ &+\sqrt{\frac{\epsilon}{\uppi}}\int_{t_1}^{t_2}
\Tr{\Gamma}{\dot{A}_t\expe{-\epsilon A^2_t}} \differ t.
\end{split}
\end{equation}  
\end{prop}

\begin{proof}
The left hand side of \eqref{gammaspectrunc1} is independent of the choice
of the unitary transformation by Proposition \ref{sf-sf-modified-formal}.  
In the right hand side of \eqref{gammaspectrunc1},  
since $U$ and $UA_t\expe{-sA_t^2}$ are bounded, 
we can rearrange these operators under the $\Upgamma$-trace in cyclic order such that
\begin{equation}\label{iso-invariance-trace}
\begin{split}
\Tr{\Upgamma}{\check{A}_t\expe{-s\check{A}^2_t}}&=\Tr{\Upgamma}{U(t_1,t) A_t\expe{-sA^2_t}U(t,t_1)}
=\Tr{\Upgamma}{A_t\expe{-s A^2_t}}.
\end{split}
\end{equation}
Hence,  in all terms,  except in $\Tr{\Gamma}{\dot{\check{A}}_t\expe{-\epsilon \check{A}^2_t}}$
we can replace $\check{A}_t$ by $A_t$.  For the final term,  differentiating $\check{A}_t$ yields
\begin{align*}
\Tr{\Gamma}{\dot{\check{A}}_t\expe{-\epsilon \check{A}^2_t}} 
&= \Tr{\Upgamma}{\dot{U}(t_1,t)A_t\expe{-sA_t^2}U(t,t_1)}
\\ &+\Tr{\Upgamma}{U(t_1,t)A_t\expe{-sA_t^2}\dot{U}(t,t_1)}
\\ &+\Tr{\Gamma}{U(t_1,t)\dot{A}_t\expe{-\epsilon A^2_t}U(t,t_1)}
\end{align*}
Differentiating $U(t_1,t)U(t,t_1)=\Iop{\mathscr{H}(t_1)}$ with respect to $t$ yields
$\dot{U}(t,t_1)Ut_1,t)=-U(t,t_1)\dot{U}(t_1,t)$. 
Since $U$ and $\dot{U}A_t\expe{-sA_t^2}$ are bounded, 
we can rearrange under the $\Upgamma$-trace in cyclic order such that the first two terms cancel each other.  
Using \eqref{iso-invariance-trace} we conclude 
\begin{align*}
\Tr{\Gamma}{\dot{\check{A}}_t\expe{-\epsilon \check{A}^2_t}} 
= \Tr{\Gamma}{\dot{A}_t\expe{-\epsilon A^2_t}}.
\end{align*}
This proves the claim.
\end{proof}

We wish to relate the integrand in \eqref{gammaspectrunc1} 
with the truncated $\Upgamma$-eta invariant 
$\eta^{<\epsilon}_{\Upgamma}(A_t)$ and the eta invariant $\eta(\underline{A}_t)$ of the
corresponding Dirac operator $\underline{A}_t$ on the base.
\begin{lem}\label{derivgammaeta}
Away from values of $t$ where eigenvalues of $\underline{A}_t$ cross the origin
\begin{equation*}
\frac{\differ}{\differ t} \eta^{<\epsilon}_{\Upgamma}(A_t)=
2 \, \sqrt{\frac{\epsilon}{\uppi}}\Tr{\Upgamma}{\dot{A}_t\expe{-\epsilon A^2_t}}\,+
\frac{\differ}{\differ t} \eta(\underline{A}_t).
\end{equation*}
Irrespective the eigenvalues crossing the origin, we have for all\footnote{As noted in
Proposition \ref{propspecflow}, the derivative $\frac{\differ}{\differ t} \xi(\underline{A}_t)$ is smooth,
while $\xi(\underline{A}_t)$ is not.} $t \in [t_1,t_2]$
\begin{equation*}
\frac{\differ}{\differ t} \eta^{<\epsilon}_{\Upgamma}(A_t)=
2 \, \sqrt{\frac{\epsilon}{\uppi}}\Tr{\Upgamma}{\dot{A}_t\expe{-\epsilon A^2_t}}\,+
2 \frac{\differ}{\differ t} \xi(\underline{A}_t).
\end{equation*}
\end{lem}

\begin{proof} 
We first show a similar relation for the upper truncated $\Upgamma$-eta 
function $\eta^{<\epsilon}_{\Upgamma}(z;A_t)$.  We will use the auxilliary evolved operator 
\begin{equation}\label{A-evolved-natural}
\check{A}_t:=U(t_1,t)A_t\,U(t,t_1), \quad 
\expe{-s\check{A}^2_t}=U(t_1,t)\expe{-sA^2_t}U(t,t_1),
\end{equation}
where as before,  $U(t,t_1): L^2(\Sigma_1, \mathscr{S}_{L,E}(\Sigma_1))
\to L^2(\Sigma_2, \mathscr{S}_{L,E}(\Sigma_2))$ is the natural
isometry of Hilbert $\Upgamma$-modules,  defined by parallel transport  
Note that the isometry ensures that the operator $\check{A}_t$ 
acts on a fixed domain.  Our starting point are the equalities (cf.  \eqref{iso-invariance-trace})
\begin{equation}\label{tr-diff}
\begin{split}
&{\partial_t}\Tr{\Upgamma}{A_t\expe{-sA^2_t}}
\equiv {\partial_t}\Tr{\Upgamma}{\check{A}_t\expe{-s\check{A}^2_t}}
=\Tr{\Upgamma}{{\partial_t}(\check{A}_t\expe{-s\check{A}^2_t})}, \\
&\partial_s \left[s^{\frac{z\pm 1}{2}}\Tr{\Upgamma}{\dot{A}_t\expe{-sA^2_t}}\right]
=\Tr{\Upgamma}{{\partial_s}(s^{\frac{z\pm 1}{2}}\dot{A}_t\expe{-sA^2_t})}. 
\end{split}
\end{equation}
Note that in the first equality, the operators $\check{A}_t$ have fixed domain in $t$ and we differentiate in 
$t$. In the second equality, we have operators with domains fixed in $s$ (but not in $t$) and we differentiate in $s$. 
The equalities hold, since the Schwartz kernels of $\check{A}_t\expe{-s\check{A}^2_t}$
and $\dot{A}_t\expe{-sA^2_t}$, as well as their derivatives in $t$ and $s$, respectively, are $\Upgamma$-Hilbert Schmidt and 
in particular $\Upgamma$-trace class. We can now compute
\begin{eqnarray*}
\partial_t\Tr{\Upgamma}{A_t\expe{-sA_t^2}}&\stackrel{\eqref{tr-diff}}{=}&\Tr{\Upgamma}{\partial_t\left[\check{A}_t\expe{-s\check{A}_t^2}\right]}=\Tr{\Upgamma}{\partial_t\left[U(t_1,t)A_t\expe{-sA_t^2}U(t,t_1)\right]} \\
&=&\Tr{\Upgamma}{\dot{U}(t_1,t)A_t\expe{-sA_t^2}U(t,t_1)+U(t_1,t)A_t\expe{-sA_t^2}\dot{U}(t,t_1)}\\
& +&\Tr{\Upgamma}{U(t_1,t)\partial_t\left[A_t\expe{-sA_t^2}\right]U(t,t_1)}.
\end{eqnarray*}
Arguing exactly as in the proof of Proposition \ref{gammaspectrunc1-prop}, 
the first two terms cancel,  the isometries in the third term cancel as well 
and we obtain
\begin{eqnarray*}
\partial_t\Tr{\Upgamma}{A_t\expe{-sA_t^2}}
&=& \Tr{\Upgamma}{\dot{A}_t\expe{-sA_t^2}+A_t{\partial_t}\expe{-sA_t^2}}\\
&=&\Tr{\Upgamma}{\dot{A}_t\expe{-sA_t^2}-2sA_t\dot{A}_tA_t\expe{-sA_t^2}}\\
&\stackrel{(*)}{=}& \Tr{\Upgamma}{\dot{A}_t\left[\expe{-sA_t^2}-2sA_t^2\expe{-sA_t^2}\right]}=\Tr{\Upgamma}{\dot{A}_t\left(1+2s{\partial_s}\right)\expe{-sA_t^2}}.
\end{eqnarray*}  
For ($\ast$) we used that 
$\expe{-sA_t^2}$ and $A_t^k\expe{-sA_t^2}$ are bounded. This allows us to commute cyclically 
under the trace.  From here we conclude, using integration by parts in the last step
\begin{eqnarray*}
&& \int_{0}^\epsilon s^{\frac{z-1}{2}} \partial_t\Tr{\Upgamma}{A_t\expe{-sA_t^2}} \differ s = \int_{0}^\epsilon s^{\frac{z-1}{2}} \Tr{\Upgamma}{\dot{A}_t\left(1+2s{\partial_s}\right)\expe{-sA_t^2}} \differ s \\
&=& \int_{0}^\epsilon s^{\frac{z-1}{2}} \Tr{\Upgamma}{\dot{A}_t \expe{-sA_t^2}} \differ s + 2\int_{0}^\epsilon s^{\frac{z+1}{2}} {\partial_s} \Tr{\Upgamma}{\dot{A}_t\expe{-sA_t^2}} \differ s \\
&=&\int_{0}^\epsilon s^{\frac{z-1}{2}} \Tr{\Upgamma}{\dot{A}_t \expe{-sA_t^2}} \differ s + 2 s^{\frac{z+1}{2}}\Tr{\Upgamma}{\dot{A}_t\expe{-sA_t^2}}\Bigr\vert_{0}^\epsilon \\
&&\quad - (z+1)\int_{0}^\epsilon s^{\frac{z-1}{2}} \Tr{\Upgamma}{\dot{A}_t\expe{-sA_t^2}} \differ s.
\end{eqnarray*}
Using asymptotic expansion of $\Tr{\Upgamma}{\dot{A}_t\expe{-sA_t^2}}$ as $s \to 0$
we conclude for $\Rep{z}$ sufficiently large
\begin{equation}\label{side1}
\begin{split}
\int_{0}^\epsilon s^{\frac{z-1}{2}} \partial_t\Tr{\Upgamma}{\check{A}_t\expe{-s\check{A}_t^2}} \differ s 
&= 2 \epsilon^{\frac{z+1}{2}}\Tr{\Upgamma}{\dot{A}_t\expe{-\epsilon A_t^2}}
\\ &-z \int_{0}^\epsilon s^{\frac{z-1}{2}} \Tr{\Upgamma}{\dot{A}_t\expe{-sA_t^2}} \differ s.
\end{split}
\end{equation}
On the other hand, we conclude with the dominated convergence theorem
\begin{equation*}
\begin{split}
\int_{0}^\epsilon s^{\frac{z-1}{2}} \partial_t\Tr{\Upgamma}{A_t\expe{-sA_t^2}} \differ s 
&= \frac{\differ}{\differ t} \int_{0}^\epsilon s^{\frac{z-1}{2}}
\Tr{\Upgamma}{A_t\expe{-sA_t^2}} \differ s 
\\ &= \Gamma\left(\frac{z+1}{2}\right)\frac{\differ}{\differ t}\eta^{< \epsilon}_\Upgamma(z; A_t).
\end{split}
\end{equation*}
Comparing this with \eqref{side1}, we obtain
\begin{equation}\label{etaetabasefuncz1}
\frac{\differ}{\differ t}\eta^{< \epsilon}_\Upgamma(z;A_t)= \frac{2\epsilon^{\frac{z+1}{2}}}{\Gamma\left(\frac{z+1}{2}\right)}\Tr{\Upgamma}{\dot{A}_t\expe{-\epsilon A_t^2}}-\frac{z}{\Gamma\left(\frac{z+1}{2}\right)}\int_{0}^\epsilon s^{\frac{z-1}{2}} \Tr{\Upgamma}{\dot{A}_t\expe{-sA_t^2}} \differ s. 
\end{equation}
In order to study the second summand, we note that the pointwise asymptotics of 
the Schwartz kernels of $\dot{A}_t \expe{-sA_t^2}$ and 
on the cover $\Sigma$ and of $\dot{\underline{A}}_t\expe{-s\underline{A}_t^2}$ 
on the base coincide. Namely for any $p \in \mathscr{F} \cong \Sigma_\Gamma$, 
where the identification of the fundamental domain $\mathscr{F}$ with the base 
$\Sigma_\Gamma$ holds away from a set of measure zero, we have for the pointwise traces
\begin{equation}\label{asymptdotAkernelerr}
\begin{split}
&\tr{p}{\dot{A}_t \expe{-sA_t^2}(p,p,t;s)}
\quad\AAsympt{s} \quad
\sum_{j\in \N_0} b_j(p,t) s^{\frac{- n+j-1}{2}}, \\ 
&\tr{p}{\dot{\underline{A}}_t\expe{-s\underline{A}_t^2}(p,p,t;s)}\quad\AAsympt{s} \quad
\sum_{j\in \N_0} b_j(p,t) s^{\frac{- n+j-1}{2}}.
\end{split}
\end{equation}
Consequently, we obtain
\begin{equation*}
\begin{split}
\underset{z=0}{\textup{Res}}\int_{0}^\epsilon s^{\frac{z-1}{2}} \Tr{\Upgamma}{\dot{A}_t\expe{-sA_t^2}} \differ s 
= 2 \int_{\mathscr{F}} b_n(p,t) = \underset{z=0}{\textup{Res}}\int_{0}^\infty s^{\frac{z-1}{2}} 
\Tr{}{\dot{\underline{A}}_t\expe{-s\underline{A}_t^2}} \differ s.
\end{split}
\end{equation*}
We conclude in view of \eqref{etaetabasefuncz1}
\begin{equation}\label{etaetabasefuncz2}
\frac{\differ}{\differ t}\eta^{< \epsilon}_\Upgamma(A_t)=
2\sqrt{\frac{\epsilon}{\uppi}}\Tr{\Upgamma}{\dot{A}_t\expe{-\epsilon A_t^2}}-
\frac{1}{\sqrt{\uppi}} \, \underset{z=0}{\textup{Res}}\int_{0}^\infty s^{\frac{z-1}{2}} 
\Tr{}{\dot{\underline{A}}_t\expe{-s\underline{A}_t^2}} \differ s
\end{equation}
The statement now follows from
Proposition 1.10.3 (b) in \cite{gilkeyinvariance}.
\end{proof}

\begin{cor}\label{sf-sf-modified} 
Consider the smooth family of twisted hypersurface 
Dirac operators $A_t := A^{E_L}_t$ for $t \in I =[t_1,t_2]$.  Then the 
$\Gamma$-spectral flow in Definition \ref{defimodgammaspectralflow}
satisfies the following expression
\begin{eqnarray}\label{gammaspecflowwithxi}
\specflow{A_{t}}{t \in [t_1,t_2]} = 
 \xi_\Upgamma(A_{2})-\xi_\Upgamma(A_{1})-\int_{t_1}^{t_2}\frac{\differ}{\differ t} 
 \xi(\underline{A}_t)\differ t.
\end{eqnarray}
\end{cor}

\begin{proof}
Lemma \ref{derivgammaeta} yields
\begin{equation*}
\begin{split}
\sqrt{\frac{\epsilon}{\uppi}}\int_{t_1}^{t_2}\Tr{\Upgamma}{\dot{A}_t\expe{-\epsilon A^2_t}}\differ t 
&=\frac{1}{2}\int_{t_1}^{t_2}\frac{\differ}{\differ t} \eta^{<\epsilon}_{\Upgamma}(A_t)\differ t
-\int_{t_1}^{t_2}\frac{\differ}{\differ t} \xi(\underline{A}_t)\differ t \\
&=\frac{\eta^{<\epsilon}_\Upgamma(A_{2})-\eta^{<\epsilon}_\Upgamma(A_{1})}{2}
-\int_{t_1}^{t_2}\frac{\differ}{\differ t} \xi(\underline{A}_t)\differ t.
\end{split}
\end{equation*}
Comparing this with \clef{gammaspectrunc1} shows
\begin{eqnarray}
\specflow{A_{t}}{t \in [t_1,t_2]} = 
 \xi_\Upgamma(A_{2})-\xi_\Upgamma(A_{1})-\int_{t_1}^{t_2}\frac{\differ}{\differ t} 
 \xi(\underline{A}_t)\differ t.
\end{eqnarray}
\end{proof}

Note that the expression in \clef{gammaspecflowwithxi} is the wanted $\Upgamma$-version of Proposition \ref{propspecflow}.

\section{$L^2$-Gamma index of the Lorentzian Dirac operator}\label{chap:algind}

As above, we shall write for the smooth family of twisted hypersurface 
Dirac operators $A_t := A^{E_L}_t, t \in I =[t_1,t_2]$. 
Let us recall the central objects of this section and their notation. 
We begin with the spectral projections for $A_t$.
\begin{equation}\label{P-overline-def}
\begin{split}
P_{>0}(t) &= \Chi_{\intervallo{0}{\infty}{}}(A_t), P_{\geq 0}(t) = \Chi_{\intervallro{0}{\infty}{}}(A_t), \\
P_{<0}(t) &= \Chi_{\intervallo{-\infty}{0}{}}(A_t), P_{\leq 0}(t) = \Chi_{\intervalllo{-\infty}{0}{}}(A_t), \\
\overline{P}^{>0}_{\geq 0}(t) &:=  
\left. P_{\geq 0}(t) \right|_{\, \range{P_{>0}(t)}}:  \range{P_{>0}(t)} \to \range{P_{\geq 0}(t)}.
\end{split}
\end{equation}  
The spectral projections for the evolved 
operators $\hat{A}_t:=Q(t_1,t)A_t\,Q(t,t_1)$ are denoted 
as in the previous section by an additional upper hat.
\begin{equation}
\begin{split}
\hat{P}_{>0}(t) &= \Chi_{\intervallo{0}{\infty}{}}(\hat{A}_t), \hat{P}_{\geq 0}(t) = \Chi_{\intervallro{0}{\infty}{}}(\hat{A}_t), \\
\hat{P}_{<0}(t) &= \Chi_{\intervallo{-\infty}{0}{}}(\hat{A}_t), \hat{P}_{\leq 0}(t) = \Chi_{\intervalllo{-\infty}{0}{}}(\hat{A}_t), \\
\hat{\overline{P}}^{\,>0}_{\geq 0}(t) &:=  \left. \hat{P}_{\geq 0}(t)\right|_{\, \range{\hat{P}_{>0}(t)}}:  \range{\hat{P}_{>0}(t)} \to \range{\hat{P}_{\geq 0}(t)}. \end{split}
\end{equation}  
A central role in our discussion here is played by the $\Gamma$-Fredholm operators 
which were introduced in Theorem \ref{qgengammafred}.
\begin{equation}
\begin{split}
Q_{++}(t_2,t_1) &:= P_{>0}(t_2)\circ Q(t_2,t_1)\circ P_{\geq 0}(t_1):
\range{P_{\geq 0}(t_1)} \to \range{P_{> 0}(t_2)}, \\
Q_{--}(t_2,t_1) &:= P_{\leq 0}(t_2)\circ Q(t_2,t_1)\circ P_{< 0}(t_1):
\range{P_{< 0}(t_1)} \to \range{P_{\leq 0}(t_2)}.
\end{split}
\end{equation} 
The main result of this section is the following theorem. 
\begin{theo}\label{propgammaindexspecflow}
\begin{equation}\label{Qindices}
\Index_\Gamma\left(Q_{\pm\pm}(t_2,t_1)\right)=
\mp\specflow{A_t}{t\in[t_1,t_2]}\pm\dim_\Gamma\kernel{A_{2}}.
\end{equation}
\end{theo}  
The remainder of the section is devoted to the proof of this result, 
which we attain after a sequence of intermediate auxiliary lemmas.

\begin{lem}\label{evolvprojgammafred}
For any pair $\tau, t\in [t_1,t_2]$ the operators
\begin{equation*}
\begin{split}
\hat{P}_{>0}(t) \equiv \left. \hat{P}_{>0}(t) \right|_{\range{\hat{P}_{\geq 0}(\tau)}}&: 
\,\,\range{\hat{P}_{\geq 0}(\tau)}\quad \rightarrow \quad \range{\hat{P}_{>0}(t)} \\
\hat{P}_{\leq 0}(t) \equiv \left. \hat{P}_{\leq 0}(t) \right|_{\range{\hat{P}_{< 0}(\tau)}}&:
\,\,\range{\hat{P}_{< 0}(\tau)}\quad \rightarrow \quad \range{\hat{P}_{\leq 0}(t)}
\end{split}
\end{equation*}
are $\Gamma$-Fredholm with $\Gamma$-indices
\begin{equation}\label{indexevolvedproj}
\begin{split}
\Index_\Gamma\left(\left. \hat{P}_{>0}(t) \right|_{\range{\hat{P}_{\geq 0}(\tau)}}\right)
&=\Index_\Gamma\left(Q_{++}(t,\tau)\right),  \\
\Index_\Gamma\left(\left. \hat{P}_{\leq 0}(t) \right|_{\range{\hat{P}_{< 0}(\tau)}}\right)
&=\Index_\Gamma\left(Q_{--}(t,\tau)\right).
\end{split}
\end{equation}
In particular,  we recover
\begin{equation}\label{Qplusminus}
\Index_\Gamma\left(Q_{++}(t,\tau)\right)
= - \Index_\Gamma\left(Q_{--}(t,\tau)\right)
\end{equation}
from \cite[Theorem 7.5]{OD2}.
\end{lem} 
\begin{proof}
Let $v\in\range{P_{\geq 0}(\tau)}$ and $u=Q(t_1,\tau)v \in \range{\hat{P}_{\geq 0}(\tau)}$, then
\begin{equation*}
\begin{split}
\hat{P}_{>0}(t)u&=Q(t_1,t)Q_{++}(t,\tau)v=Q(t_1,t)Q_{++}(t,\tau)Q(\tau,t_1)u.
\end{split}
\end{equation*}
Similarly, if $v\in\range{P_{< 0}(\tau)}$ and $u=Q(t_1,\tau)v \in \range{\hat{P}_{< 0}(\tau)}$, we get
\begin{equation*}
\begin{split}
\hat{P}_{\leq 0}(t)u&= Q(t_1,t)Q_{--}(t,\tau)v=Q(t_1,t)Q_{--}(t,\tau)Q(\tau,t_1)u.
\end{split}
\end{equation*}
Since $Q_{\pm\pm}$ are $\Gamma$-Fredholm for all subintervals of $[t_1,t_2]$,
and $Q$ are unitary for all time intervals with vanishing $\Gamma$-index, 
the statement \eqref{indexevolvedproj} follows. For the final observation, we compute
using Proposition \ref{propindexfredpair} (3) in the second equality, and 
Lemma \ref{evolvprojgammafred} in the final equality
\begin{equation*}
\begin{split}
\Index_\Gamma\left(Q_{++}(t,\tau)\right)
&= \Index_\Gamma\left(\left. \hat{P}_{>0}(t) \right|_{\range{\hat{P}_{\geq 0}(\tau)}}\right)
= - \Index_\Gamma\left(\left. \mathds{1} - \hat{P}_{>0}(t) \right|_{\range{\mathds{1} - \hat{P}_{\geq 0}(\tau)}}\right)
\\ &= - \Index_\Gamma\left(\left. \hat{P}_{\leq 0}(t) \right|_{\range{\hat{P}_{< 0}(\tau)}}\right)
= - \Index_\Gamma\left(Q_{--}(t,\tau)\right).
\end{split}
\end{equation*}
\end{proof}

\begin{cor}
For any pair $\tau, t\in [t_1,t_2]$ the operators
\begin{equation*}
\begin{split}
\hat{P}_{\geq 0}(t) \equiv \left. \hat{P}_{\geq 0}(t) \right|_{\range{\hat{P}_{\geq 0}(\tau)}}&: 
\,\,\range{\hat{P}_{\geq 0}(\tau)}\quad \rightarrow \quad \range{\hat{P}_{\geq 0}(t)} \\
\hat{P}_{< 0}(t) \equiv \left. \hat{P}_{< 0}(t) \right|_{\range{\hat{P}_{< 0}(\tau)}}&:
\,\,\range{\hat{P}_{< 0}(\tau)}\quad \rightarrow \quad \range{\hat{P}_{< 0}(t)}
\end{split}
\end{equation*}
are $\Gamma$-Fredholm with $\Gamma$-indices 
\begin{equation}\label{indexevolvedproj2}
\begin{split}
\Index_\Gamma\left(\left. \hat{P}_{\geq 0}(t) \right|_{\range{\hat{P}_{\geq 0}(\tau)}}\right)
&=\Index_\Gamma\left(Q_{++}(t,\tau)\right)-\dim_\Gamma\kernel{A_t},  \\
\Index_\Gamma\left(\left. \hat{P}_{< 0}(t) \right|_{\range{\hat{P}_{< 0}(\tau)}}\right)
&=  \Index_\Gamma\left(Q_{--}(t,\tau)\right) + \dim_\Gamma\kernel{A_t}.
\end{split}
\end{equation}
\end{cor}

\begin{proof}
Consider the projection $P_{\{0\}}(t)$ onto $\ker A_t$. We compute
\begin{align*}
P_{\geq 0}(t) &= P_{\geq 0}(t)^2 = P_{\geq 0}(t) \left(P_{> 0}(t) + P_{\{0\}}(t)\right) 
\\ &= P_{\geq 0}(t) P_{> 0}(t) + P_{\geq 0}(t)  P_{\{0\}}(t)
\\ &= \overline{P}^{>0}_{\geq 0}(t) P_{> 0}(t) + P_{\{0\}}(t).
\end{align*}
Similarly, we obtain for the evolved operators
\begin{align*}
\left. \hat{P}_{\geq 0}(t) \right|_{\range{\hat{P}_{\geq 0}(\tau)}}
= \left. \hat{\overline{P}}^{>0}_{\geq 0}(t) \hat{P}_{> 0}(t)\right|_{\range{\hat{P}_{\geq 0}(\tau)}} + 
\left. \hat{P}_{\{0\}}(t)\right|_{\range{\hat{P}_{\geq 0}(\tau)}}.
\end{align*}
The restricted projector $\overline{P}^{\,>0}_{\geq 0}(t)$ is $\Upgamma$-Fredholm for all $t$ which can be directly proven with the help of \cite[Lem. 7.1]{OD2} or \cite[Lem.8.16]{ODT}.
Since $Q$ ist unitary for all times, the evolved restricted projector $\hat{\overline{P}}^{>0}_{\geq 0}(t)$ becomes $\Gamma$-Fredholm, too. The previous Lemma then implies that the composition on the righthand-side is $\Gamma$-Fredholm.
The projection $P_{\{0\}}(t)$ onto $\ker A_t$ is $\Gamma$-trace class. Indeed,
$P_{\{0\}}(t) = P_{\{0\}}(t) \circ e^{-A_t^2}$, where $P_{\{0\}}(t)$ is bounded and $e^{-A_t^2}$
is trace class. Hence the evolved projection $\hat{P}_{\{0\}}(t)$ is trace-class as well and we
conclude by Lemma \ref{evolvprojgammafred}
\begin{align*}
\Index_\Gamma\left( \left. \hat{P}_{\geq 0}(t) \right|_{\range{\hat{P}_{\geq 0}(\tau)}} \right)
&= \Index_\Gamma\left( \left. \hat{\overline{P}}^{>0}_{\geq 0}(t) \hat{P}_{> 0}(t)\right|_{\range{\hat{P}_{\geq 0}(\tau)}}\right)
\\ &= \Index_\Gamma\left( \hat{\overline{P}}^{>0}_{\geq 0}(t)\right)
+ \Index_\Gamma\left( \left. \hat{P}_{> 0}(t)\right|_{\range{\hat{P}_{\geq 0}(\tau)}}\right)
\\ &= \Index_\Gamma\left( \hat{\overline{P}}^{>0}_{\geq 0}(t)\right) + 
\Index_\Gamma\left(Q_{++}(t,\tau)\right)
\\ &= \Index_\Gamma\left( \overline{P}^{>0}_{\geq 0}(t)\right) + 
\Index_\Gamma\left(Q_{++}(t,\tau)\right).
\end{align*}
The last equivalence is justified since $Q$ is unitary for all time intervals with vanishing $\Gamma$-index. We can identify the $\Gamma$-index
of the restricted projections: by definition, we have for any $t \in [t_1,t_2]$ (note that
$A_t$ is $\Gamma$-Fredholm and hence its kernel $\ker A_t$ is of 
finite $\Gamma$-dimension)
\begin{equation}\label{gamindrestprojdimkern}
\begin{split}
&\Index_\Gamma\left(\overline{P}^{\geq 0}_{> 0}(t)\right)=\Index_\Gamma\left(\left. P_{>0}(t) \right|_{\range{P_{\geq 0}(t)}}\right)
= \dim_\Gamma \ker A_{t}, \\
&\Index_\Gamma\left(\overline{P}^{> 0}_{\geq 0}(t)\right)=\Index_\Gamma\left(\left. P_{\geq 0}(t) \right|_{\, \range{P_{>0}(t)}}\right)
= -\dim_\Gamma \ker A_{t}.
\end{split}
\end{equation}
The first statement in \eqref{indexevolvedproj2} then follows. The second statement follows
by Proposition \ref{propindexfredpair} (3) and \eqref{Qplusminus}.
\end{proof}

We can now prove the main result of this section.

\begin{proof} (Theorem \ref{propgammaindexspecflow}) 
We continue in the notation of Proposition \ref{sf-sf-modified}.
Using Proposition \ref{propindexfredpair} (2) and (4) we compute
\begin{align*}
\Index_\Gamma \Bigl(\hat{P}_{\geq 0}(\tau_{j-1})\Big\vert \hat{P}_{\geq 0}(\tau_{j})\Bigr)
&= \Index_\Gamma \Bigl(\hat{P}_{\geq 0}(\tau_{j-1})\Big\vert \hat{P}_{\geq 0}(t_1)\Bigr)
+ \Index_\Gamma \Bigl(\hat{P}_{\geq 0}(t_1)\Big\vert \hat{P}_{\geq 0}(\tau_{j})\Bigr)
\\ &= \Index_\Gamma \Bigl(\hat{P}_{\geq 0}(\tau_{j-1})\Big\vert \hat{P}_{\geq 0}(t_1)\Bigr)
- \Index_\Gamma \Bigl(\hat{P}_{\geq 0}(\tau_{j})\Big\vert \hat{P}_{\geq 0}(t_1)\Bigr).
\end{align*}
Plugging this into the formula \eqref{modgammaspectralflow} for the spectral flow,  we obtain
\begin{align*}
\specflow{A_t}{t\in I} &:= 
\sum_{j=1}^N\Index_\Gamma
\Bigl(\hat{P}_{\geq 0}(\tau_{j-1})\Big\vert \hat{P}_{\geq 0}(\tau_{j})\Bigr) \\
&= \sum_{j=1}^N \Index_\Gamma \Bigl(\hat{P}_{\geq 0}(\tau_{j-1})\Big\vert \hat{P}_{\geq 0}(t_1)\Bigr)
- \Index_\Gamma \Bigl(\hat{P}_{\geq 0}(\tau_{j})\Big\vert \hat{P}_{\geq 0}(t_1)\Bigr)
\\ &= \Index_\Gamma \Bigl(\hat{P}_{\geq 0}(t_1)\Big\vert \hat{P}_{\geq 0}(t_1)\Bigr) - 
\Index_\Gamma \Bigl(\hat{P}_{\geq 0}(t_2)\Big\vert \hat{P}_{\geq 0}(t_1)\Bigr).
\end{align*}
Note that by definition of relative $\Gamma$-indices
$$
\Index_\Gamma \Bigl(\hat{P}_{\geq 0}(t)\Big\vert \hat{P}_{\geq 0}(\tau)\Bigr)
\equiv \Index_\Gamma\left( \left. \hat{P}_{\geq 0}(t) \right|_{\range{\hat{P}_{\geq 0}(\tau)}} \right).
$$
Using \eqref{indexevolvedproj2}, we conclude 
\begin{equation}\begin{split}
\specflow{A_t}{t\in I} &= - \Index_\Gamma\left(Q_{++}(t_2,t_1)\right)
+ \dim_\Gamma \ker A_{t_2} \\ &= \Index_\Gamma\left(Q_{--}(t_2,t_1)\right)
+ \dim_\Gamma \ker A_{t_2}.
\end{split}\end{equation}
\end{proof}

Using Theorem \ref{indexDgaAPStwist}, notation in \eqref{Qabbreviation} and 
the observation \eqref{Qplusminus}, our result in Theorem \ref{propgammaindexspecflow}
implies directly the following Corollary. 
\begin{cor}\label{DQsf}
\begin{align*}
&\Index_{\Gamma}(D^{E_L}_{+,\mathrm{APS}(0,0)})
= \Index_\Gamma\left(Q_{--}(t_2,t_1)\right) = 
\specflow{A_t}{t\in I} - \dim_\Gamma \ker A_{t_2}, \\
&\Index_{\Gamma}(D^{E_L}_{+,\mathrm{aAPS}(0,0)})
= \Index_\Gamma\left(Q_{++}(t_2,t_1)\right)
= - \specflow{A_t}{t\in I} + \dim_\Gamma \ker A_{t_2}.
\end{align*}
\end{cor}

\section{Proof of the $L^2$-Gamma index theorem on spacetimes}\label{chap:geoind}

We continue in the geometric setting as outlined in \S \ref{setting-subsection}.
Let $\check{D}^{E_L}$ be the 
twisted spin$^{c}$ Dirac operator of $M = [t_1,t_2] \times \Sigma$ with the (flipped) Riemannian 
metric $\check{\met} = N^2 dt^2 + \met_t$, restricted to positive spinors. 
We deform the Riemannian metric $\check{\met}$ to be ultra-static near the boundary, 
i.e.  the mean curvatures of the hypersurfaces $\Sigma_t$ vanish identically
near $t\in \SET{t_1,t_2}$ and the hypersurface metric $\met_t$ becomes $\met_1=\met_{t_1}$ for $t$ in a neighborhood of $t_1$ and $\met_2=\met_{t_2}$ for $t$ in a neighborhood of $t_2$.
Then $\nu=-\partial_t$ near the boundary hypersurfaces $\Sigma_j=\SET{t_j}\times\Sigma$ with $j=1,2$.  $\check{D}^{E_L}$ then becomes 
\begin{align}
\check{D}^{E_L} = 
- (\upbeta\otimes \Iop{E_L})(\partial_t + A^{E_L}_j).
\end{align}
We denote by $\check{D}^{E_L}_{\mathrm{APS}}$ the operator with classical APS boundary conditions.  Note that we do not impose
the ultra-static assumption on the Lorentzian metric $\met$, but only on the auxiliary Riemannian metric $\check{\met}$.  The following result follows by identifying the $\Gamma$-index of $\check{D}^{E_L}_{\mathrm{APS}}$
in terms of the spectral flow $\specflow{A^{E_L}_t}{t\in I}$.  

\begin{prop}\label{riemgammaspecflow} In the notation of Theorem \ref{mainres2}
\begin{equation}
\begin{split}
\specflow{A_{t}^{E_L}}{t \in [t_1,t_2]}
 =  \int_{M_\Gamma}\hat{\mathbbm{A}}\left(\check{M}_\Gamma\right)\wedge\expe{\mathrm{c}_1
\left(L_\Gamma\right)/2}\wedge \mathrm{ch}\left(E_\Gamma\right)
 -\xi_\Gamma(A^{E_L}_1)+\xi_\Gamma(A^{E_L}_2),
 \end{split}
\end{equation}
where $\hat{\mathbbm{A}}\left(\check{M}_\Gamma\right)$ is the 
$\check{A}$-genus of $(M_\Gamma,\check{\met}_\Gamma)$.
\end{prop}

\begin{proof}
Since $\partial_t$ is inward pointing normal vector field at $\Sigma_1$, 
and outward pointing at $\Sigma_2$, the $\Gamma$-index theorem of Ramachandran
\cite[Thm.1.1]{ramachandran} asserts that $\check{D}^{E_L}_{\mathrm{APS}}$
is $\Gamma$-Fredholm with $\Gamma$-index 
\begin{equation}\label{riemgammaindexI}
\Index_\Gamma(\check{D}^{E_L}_{\mathrm{APS}})=
\int_{M_\Gamma}\hat{\mathbbm{A}}\left(\check{M}_\Gamma\right)\wedge\expe{\mathrm{c}_1
\left(L_\Gamma\right)/2}\wedge \mathrm{ch}\left(E_\Gamma\right)
-\xi_\Gamma(A^{E_L}_1)-\xi_\Gamma(-A^{E_L}_2),
\end{equation}
where the integrand is expressed in terms of characteristic 
classes on the compact base.
Note that by definition (cf. \eqref{gammaxi})
\begin{equation}
\xi_\Upgamma(- A) =\frac{\eta_\Upgamma(- A)+\dim_\Upgamma\kernel{-A}}{2}
= - \xi_\Upgamma(A) + \dim_\Upgamma\kernel{A}.
\end{equation}
Hence the boundary contribution in \clef{riemgammaindexI} can be rewritten as
\begin{equation}\label{xiplusminus}
\begin{split}
\xi_\Gamma(A^{E_L}_1)+\xi_\Gamma(-A^{E_L}_2)
&=\xi_\Gamma(A_{1}^{E_L})-\xi_\Gamma(A_{2}^{E_L})+\dim_\Gamma \kernel{A_{2}^{E_L}}.
\end{split}
\end{equation}
The $\Gamma$-spectral flow formula \clef{gammaspecflowwithxi} 
together with \eqref{riemgammaindexI} and \eqref{xiplusminus} now implies
\begin{equation}\label{riemgammaindexII}
\begin{split}
\Index_\Gamma(\check{D}^{E_L}_{\mathrm{APS}})=&\int_{M_\Gamma}\hat{\mathbbm{A}}\left(\check{M}_\Gamma\right)\wedge\expe{\mathrm{c}_1\left(L_\Gamma\right)/2}\wedge \mathrm{ch}\left(E_\Gamma\right)+\specflow{A_{t}^{E_L}}{t \in [t_1,t_2]}\\
&+ \frac{1}{2}\int_{t_1}^{t_2} \frac{\differ}{\differ t} \xi(\underline{A}_{t}^{E_L}) \differ t 
-\dim_\Gamma \kernel{A_{2}^{E_L}},
\end{split}
\end{equation}
where $\underline{A}_{t}^{E_L}$ denotes the
boundary operators of the Dirac operator $\underline{\check{D}}^{E_L}_{\mathrm{APS}}$
on the compact base $M_\Gamma \equiv M\mirror{\setminus}\Gamma =[t_1,t_2]\times \Sigma\mirror{\setminus}\Gamma$. 
Using Proposition \ref{propspecflow} 
instead of \clef{gammaspecflowwithxi}, we obtain 
\begin{equation}\label{apsindexstep1}
\begin{split}
\Index\left(\underline{\check{D}}^{E_L}_{\mathrm{APS}}\right) 
&= \int_{M_\Gamma}\hat{\mathbbm{A}}\left(\check{M}_\Gamma\right)\wedge\expe{\mathrm{c}_1\left(L_\Gamma\right)/2}\wedge \mathrm{ch}\left(E_\Gamma\right) + \Specflow{\underline{A}_{t}^{E_L}}{t \in [t_1,t_2]} \\
&+ \frac{1}{2}\int_{t_1}^{t_2} \frac{\differ}{\differ t} \xi(\underline{A}_{t}^{E_L}) \differ t 
-\dim \kernel{\underline{A}_{2}^{E_L}}. 
\end{split}
\end{equation}
Note that in both expressions, in \clef{riemgammaindexII} and \clef{apsindexstep1}, 
the $\xi$-invariant $\xi(\underline{A}_{t}^{E_L})$, defined on the base, appears.
Comparing \clef{riemgammaindexII} and \clef{apsindexstep1}, we obtain
\begin{align*}
&\Index_\Gamma(\check{D}^{E_L}_{\mathrm{APS}})-
\specflow{A_{t}^{E_L}}{t \in [t_1,t_2]}-\frac{1}{2}\int_{t_1}^{t_2} 
\frac{\differ}{\differ t} \xi(\underline{A}_{t}^{E_L}) \differ t +\dim_\Gamma \kernel{A^{E_L}_2} \\
&= \Index\left(\underline{\check{D}}^{E_L}_{\mathrm{APS}}\right)-
\Specflow{\underline{A}_{t}^{E_L}}{t \in [t_1,t_2]}- \frac{1}{2}\int_{t_1}^{t_2} 
\frac{\differ}{\differ t} \xi(\underline{A}_{t}^{E_L}) \differ t +\dim\kernel{\underline{A}_{2}^{E_L}}\\
&= \int_{M_\Gamma}\hat{\mathbbm{A}}\left(\check{M}_\Gamma\right)\wedge\expe{\mathrm{c}_1
\left(L_\Gamma\right)/2}\wedge \mathrm{ch}\left(E_\Gamma\right).
\end{align*}
Subtracting the integral term on both sides of the first equality, yields
\begin{eqnarray*}
&& \Index_\Gamma(\check{D}^{E_L}_{\mathrm{APS}})-\specflow{A_{t}^{E_L}}{t \in [t_1,t_2]}
 +\dim_\Gamma \kernel{A_{2}^{E_L}} \\
&=&  \Index\left(\underline{\check{D}}^{E_L}_{\mathrm{APS}}\right)-
\Specflow{\underline{A}_{t}^{E_L}}{t \in [t_1,t_2]} +\dim\kernel{\underline{A}_{2}^{E_L}} \quad \in \Z.
\end{eqnarray*}
Because the right-hand side of the equality is well-known to be integer-valued, the left-hand side also becomes integer-valued. Note also from \eqref{apsindexstep1} that the right hand side is continuous in $t_1, t_2$, indeed 
\eqref{apsindexstep1} rewrites as 
\begin{equation}\label{apsindexstep2}
\begin{split}
&\Index\left(\underline{\check{D}}^{E_L}_{\mathrm{APS}}\right) 
- \Specflow{\underline{A}_{t}^{E_L}}{t \in [t_1,t_2]} 
+ \dim \kernel{\underline{A}_{2}^{E_L}}
\\ &= \int_{M_\Gamma}\hat{\mathbbm{A}}\left(\check{M}_\Gamma\right)\wedge
\expe{\mathrm{c}_1\left(L_\Gamma\right)/2}\wedge \mathrm{ch}\left(E_\Gamma\right) 
+ \frac{1}{2}\int_{t_1}^{t_2} \frac{\differ}{\differ t} \xi(\underline{A}_{t}^{E_L}) \differ t.
\end{split}
\end{equation}
The integral expressions on the right hand side of \eqref{apsindexstep2} are continuous in 
$t_1, t_2$ and converge to zero for $t_2 \to t_1$. Together with the left hand side 
of  \eqref{apsindexstep2} being integer, we conclude that it is identically zero and hence
\begin{eqnarray*}
\Index_\Gamma(\check{D}^{E_L}_{\mathrm{APS}})-\specflow{A_{t}^{E_L}}{t \in [t_1,t_2]}
 +\dim_\Gamma \kernel{A_{2}^{E_L}}  = 0.
\end{eqnarray*}
We arrive in view of \eqref{riemgammaindexI}
at the following expression of the $\Gamma$-spectral flow
\begin{equation}
\begin{split}
\specflow{A_{t}^{E_L}}{t \in [t_1,t_2]}
 &= \Index_\Gamma(\check{D}^{E_L}_{\mathrm{APS}})
 + \dim_\Gamma \kernel{A_{2}^{E_L}}
 \\ &=  \int_{M_\Gamma}\hat{\mathbbm{A}}\left(\check{M}_\Gamma\right)\wedge\expe{\mathrm{c}_1
\left(L_\Gamma\right)/2}\wedge \mathrm{ch}\left(E_\Gamma\right)
 -\xi_\Gamma(A^{E_L}_1)+\xi_\Gamma(A^{E_L}_2). 
 \end{split}
\end{equation}
This proves the claim.
\end{proof}

We can finally prove our main result, Theorem \ref{mainres2},
which involves rewriting the integral term in the formula of 
Proposition \ref{riemgammaspecflow} in terms of $\met$ instead of
$\check{\met}$.

\begin{proof} (proof of Theorem \ref{mainres2}) 
Proposition \ref{riemgammaspecflow} and Corollary \ref{DQsf} imply
\begin{align*}
\Index_{\Gamma}(D^{E_L}_{+,\mathrm{APS}})
&= \int_{M_\Gamma}\hat{\mathbbm{A}}\left(\check{M}_\Gamma\right)\wedge\expe{\mathrm{c}_1
\left(L_\Gamma\right)/2}\wedge \mathrm{ch}\left(E_\Gamma\right) \\ 
&-\xi_\Gamma(A^{E_L}_1)+\xi_\Gamma(A^{E_L}_2) - \dim_\Gamma \ker A_{t_2}. 
\end{align*}
We rewrite, as in \cite[(14)]{BaerStroh}, the integral term in terms of $\met$ instead of
$\check{\met}$ with the Stokes-Cartan theorem to
\begin{align*}
\int_{M_\Gamma}\hat{\mathbbm{A}}\left(\check{M}\right)
\wedge\expe{\mathrm{c}_1\left(L_\Gamma\right)/2}\wedge \mathrm{ch}\left(E_\Gamma\right)
&= \int_{M_\Gamma} \widehat{\mathbbm{A}}\left(M_\Gamma\right)\wedge\expe{\mathrm{c}_1\left(L_\Gamma\right)/2}\wedge \mathrm{ch}\left(E_\Gamma\right) \\
&+\int_{\bound M_\Gamma} \iota^\ast(T\widehat{\mathbbm{A}}(\check{\nabla},\nabla))\wedge\expe{\mathrm{c}_1\left(L_\Gamma\vert_{\bound M}\right)/2}\wedge \mathrm{ch}\left(E_\Gamma\vert_{\bound M}\right),
\end{align*} 
where $T\widehat{\mathbbm{A}}(\check{\nabla},\nabla)$ is the usual transgression 
term with $dT\widehat{\mathbbm{A}}(\check{\nabla},\nabla) = \hat{\mathbbm{A}}\left(\check{M}\right) - \hat{\mathbbm{A}}\left(M\right)$ and $\iota$ is the inclusion $\Sigma_\Gamma \hookrightarrow M_\Gamma$. The pullback of the transgression form is well-known to be independent 
of the auxiliary metric $\check{\met}$, an argument that is worked out in detail e.g. in \cite[Sec. 10.3.2]{ODT}. 
Hence we may write 
$T\widehat{\mathbbm{A}}(\met) \equiv \iota^\ast(T\widehat{\mathbbm{A}}(\check{\nabla},\nabla))$.
This proves the statement.
\end{proof}

\section{Concluding remarks and some open problems}\label{chap:opquest}

\setcounter{subsection}{1}

Our result Theorem \ref{mainres2} extends \cite[Thm.4.1]{BaerStroh}, \cite[Thm.7.1]{BaerStroh} and the index formulas from \cite[Sec.4.1/2]{BaerHan} to a spatial non-compact situation by replacing the compact Cauchy hypersurfaces with Galois coverings, associated to a Galois group $\Gamma$. All quantities in the known index formulas transfer to their pendent in the $\Gamma$-setting. Several possible modifications can be discussed or even proposed within or beyond our setting. \medskip

\subsubsection{Finite coverings}
We assume the special case that $\Gamma$ is a finite discrete Galois group such that the covering map becomes an ordinary $l$-fold covering with cardinality $l=\vert\Gamma\vert$. The Cauchy hypersurface $\Sigma$ as Galois covering with closed becomes itself compact such that the globally hyperbolic manifold $M$ becomes compact with boundary as its pendent on the covering base $M_\Gamma$. Any $\Gamma$-trace class operator becomes trace class in the ordinary sense with $l \Tr{\Gamma}{\bullet}=\Tr{}{\bullet}$. This scaling carries over to the $\Gamma$-dimension, the $\Gamma$-eta invariant and thus the $\Gamma$-xi invariant and the $\Gamma$-index. \\
\\
Considering ordinary APS boundary conditions the index of the Lorentzian Dirac operator $D^{E_L}_{+,\mathrm{APS}}$ on the covering and the $l$-fold index of $\underline{D}^{E_L}_{+,\mathrm{APS}}$ on the base differ in some boundary contributions:
\begin{equation*}
\Index\left(D^{E_L}_{+,\mathrm{APS}}\right)-l\,\Index\left(\underline{D}^{E_L}_{+,\mathrm{APS}}\right)=\left(\xi(-A_2^{E_L})-l\,\xi(-\underline{A}^{E_L}_2)\right)-l\,\left(\xi(A_1^{E_L})-l\,\xi(\underline{A}^{E_L}_1)\right).
\end{equation*}
This simple version of a Lorentzian finite covering index formula reduces for $l=1$ to $\Index\left(D^{E_L}_{+,\mathrm{APS}}\right)=\Index\left(\underline{D}^{E_L}_{+,\mathrm{APS}}\right)$ since the differences of xi-invariants coincide. Thus we have recovered the index formuala from \cite{BaerStroh}. The calculation has been worked out in detail and in the most general considered setting in \cite[Sec. 10.3.3]{ODT}.

\subsubsection{Other groups $\Gamma$}

In the present work, the group $\Gamma$ has been taken to be a discrete group of deck transformations with compact quotient $\Sigma/\Gamma$. One may want to replace the group $\Gamma$ with a locally compact unimodular group $\group{G}$ which then acts properly, cocompactly and isometrically on the complete hypersurfaces $\Sigma$. 
\medskip

The corresponding $L^2$-index theory is presented in e.g. \cite[Sec.2/3]{wang} and is quite similar to the one we used for our $\Gamma$-setting (see \cite[Sec.3/4]{OD2}) apart from a cutoff-function along the orbits which comes with a proper and cocompact $\group{G}$-action and contributes to the $\group{G}$-index. Because of this resemblance, one could have the impression that the Fredholm part in our proof can be extended to this situation which is in fact clear up to the point where the extension of Seeley's theorem on complex powers to Galois coverings is applied. 
\medskip

Our $\Gamma$-spectral flow concept can be carried over as we introduced it as a special case of spectral flow in a general semifinite von Neumann algebra. A lower truncated eta-invariant $\eta^{>\epsilon}_{\group{G}}$ can be defined as in \cite[Prop.3.3]{AzWahl} with the $\group{G}$-trace, introduced in \cite[Sec.3.1]{wang}. It is left to check that the limit $\epsilon \rightarrow 0^{+}$ and thus the full $\group{G}$-eta invariant is well-defined. The spectral flow expression similar to \clef{gammaspectrunc1} should then follow. \medskip

In order to express the ${\group{G}}$-index with geometric data, a corresponding ${\group{G}}$-index for the Riemannian Dirac operator is needed.  This has been so far worked in \cite{HoWaWa1} and \cite{HoWaWa2} out for assumption that the boundary Dirac operator is not just Fredholm, but invertible.
\medskip

Instead of replacing the group, one could also fix an element $\mathsf{g}\in\Gamma\setminus\SET{\upepsilon}$, which have finitely many conjugates in $\Gamma$, and take the finite conjugacy class $\langle \mathsf{g}\rangle$ instead of $\Gamma$. In other words, we relax the condition that the discrete group is an i.c.c. group. The case for non-compact manifolds has been studied in \cite{Lueck}. The extension to manifolds with boundary could be worked out as in \cite{ramachandran} and the discussed $\Gamma$-eta invariant reduces to Lott's delocalised eta-invariant for finite $\langle\mathsf{g}\rangle$ or of polynomial growth if $\absval{\langle\mathsf{g}\rangle}=\infty$, see \cite{LOTT19991}.
\medskip

\subsubsection{Alternative argument using Mishchenko-Fomenko bundle}

We want to point out that there is an alternative way in proving $L^2$-index theorems with respect to a Galois covering $M\rightarrow M_\Gamma$. Instead of working with the lifted Dirac operators on the covering, one can consider the operator on the compact base $M_\Gamma$, but one needs to twist the operator with the \textit{Mishchenko-Fomenko bundle} of the (reduced) $C^\ast$-algebras of $\Gamma$. 
\medskip

This treatment is suitable for a $K$-theoretic approach and proof of an index theorem as one can work on compact manifolds. On the other hand it provides a base for \textit{higher index theories} coming with topological and geometric applications like the question of the existence of positive-scalar-curvature metrics on $M$.  An example of such a higher index theorem is the Connes-Moscovici index theorem which extends the $\Gamma$-index theorem for the boundaryless case to group cochains on the group $\Gamma$; see \cite{CoMo} for the technical and conceptual details. The case with boundary, which might be interesting for the here proven $\Gamma$-index has been treated in \cite{Fu}. The spectral flow as key ingredient in our proof has a pendent for such higher index theorems, introduced in \cite{daixhang}.

\subsubsection{Further extensions}
Several conceivable extensions of the Lorentzian index theorem has been proposed in \cite{BaerStroh}. So far one has considered several further non-compact versions: the first spatial non-compact result has been presented by by Bravermann in \cite{Braverman2020} for Callias-Dirac operators on odd-dimensional globally hyperbolic manfolds. A local index theorem has been again proven by Bär and Strohmaier in \cite{BaerStroh2} for non-compact Cauchy hypersurfaces. Another non-compactness result is presented in \cite{shenwrochna} wherein the temporal compactness is relaxed and hypersurfaces at temporal infinity can be taken into account under the assumption that there the metric decays to an ultra-static metric.
\medskip

The next obvious possible extension is given by non-compact Cauchy hypersurfaces, assumed to be submanifolds of bounded geometry. Our treatment, especially in \cite{OD2}, already relied on some facts and conceptual results for this more general spatial non-compact situation, e.g. the treatment of Galois coverings as special cases of manifolds with bounded geometry in order to specify the bounded geometry version of Seeleys theorem of complex powers to our $\Gamma$-setting. The well-posedness results and the existence of the wave evolution operators has been already shown for the general case of complete Cauchy hypersurfaces in \cite{OD1}. Another spatial non-compact situation, which can be treated within a von Neumann algebra framework, is to consider the Cauchy hypersurfaces as compact foliated manifolds with non-compact leaves. Such a longitudinal foliation index has been proven by Connes and Skandalis in \cite{CoSka} for the boundaryless case and contains the family index theorem of Atiyah and Singer as special case when one replaces the leaves with fibres of a foliation. The case with boundary has been proven by Ramachandran parallely to the $\Gamma$-index case in \cite{ramachandran}. So it is naturally to ask whether and how one can adapt this situations to the globally hyperbolic case.
\medskip

A less interesting, but conceivable modification is an odd-dimensional index theorem for spacetimes as Lorentzian pendent to the thoughts, presented in \cite{Freed}. The extensions to Galois coverings and the connection of a spectral flow in this setting and the odd covering index have been investigated by Hurder in \cite{hurder}.

\end{document}